\documentclass[12pt,a4paper]{article}
\usepackage[T2A]{fontenc}
\usepackage[utf8]{inputenc}
\usepackage[russian, english]{babel}
\usepackage{epsfig,amssymb,color,graphics}
\usepackage{ifpdf}
\ifpdf
  \usepackage{epstopdf}
\fi
\usepackage{amsthm, amsmath, amssymb, amsbsy, amsfonts, latexsym, euscript}
\usepackage{enumitem}
\usepackage{authblk}
\usepackage{tikz}
\usepackage{caption}
\usetikzlibrary{arrows, topaths}
cm \hoffset = -1.0 true cm

\newtheorem{theo}{Theorem}
\newtheorem{stat}{Statement}
\newtheorem{cor}{Corollary}
\newtheorem{prp}{Proposition}
\newtheorem{lem}{Lemma}

\newtheorem{rem}{Remark}
\newenvironment{proofn}{{\bf Proof: }}
{\hfill $\diamond$\medskip}
\newcommand{\Int}{{\rm int}\,}
\newcommand{\Cl}{{\rm cl}\,}
\newcommand{\Dim}{{\rm dim}\,}
\newcommand{\Id}{{\rm id}}
\begin{document}

\title{\Large{On  Topology of Carrying Manifolds of Regular Homeomorphisms}}
\author{Elena~Gurevich\footnote{Email address: egurevich@hse.ru (E.~Gurevich)}}

\author{Ilya~Saraev\footnote{Corresponding author, Email address: isaraev@hse.ru (I.~Saraev)}}
\affil{NRU HSE, Nizhny Novgorod, Russian Federation}
\date{}

\maketitle
\abstract{We describe interrelations between a topology structure of closed manifolds (orientable and non-orientable) of the dimension $n\geq 4$ and the structure of the  non-wandering set of  regular  homeomorphisms, in particular,   Morse-Smale diffeomorphisms.}



\tableofcontents
\newpage

\section*{Introduction and statement of results}\label{intro}
\addcontentsline{toc}{section}{ Introduction and statement of results}

In 1937, A.~Andronov and L.~Pontryagin introduced the notion of  roughness  of  dynamical systems and showed that necessary and sufficient conditions of the  roughness of a flow on the 2-dimensional sphere are finiteness and hyperbolicity of its non-wandering set and the
absence of trajectories joining two saddle equilibria.  In 1960, S.~Smale introduced a similar class of dynamical systems on closed manifolds of an arbitrary dimension, for which the condition of the  absence of heteroclinic trajectories was replaced by the condition of the   transversality  of the  intersection of  invariant manifolds. In~\cite{Els, S2} for such systems the inequalities   connecting the number of periodic orbits of different types  and Betty numbers of the carrying manifold were obtained, similar to Morse inequalities for Morse functions. Since then  the systems   got a  name   \emph{Morse-Smale}.  

In particular, in~\cite{S2}, the following generalization of Poincare-Hopf formula was obtained. Let $M^n$ be a closed connected smooth manifold with  Euler characteristic  $\chi(M^n)$, $f: M^n\to M^n$ be a Morse-Smale diffeomorphism,  and  $k_i$  denotes the number of periodic points of $f$  whose unstable manifolds have dimension  $i\in \{0,\dots, n\}$.   Then   
\begin{equation}\sum \limits_{i=0}^n(-1)^ik_i=\chi(M^n).\end{equation}

 Since   Euler characteristic is a complete topological  invariant for orientable and non-orientable two-dimensional closed manifolds, the formula  above gives remarkable interrelation between the structure of the non-wandering set of a Morse-Smale system and its ambient manifold of dimension two. In particular, a genus $g$ of an orientable closed surface $M^2$ can be expressed as  \begin{equation}g=(k_1-k_0-k_2+2)/2.\end{equation} For $n\geq 3$, this  formula is  not so informative, in particular  because  $\chi(M^{2k+1})=0$ for any manifold $M^{2k+1}$ of odd dimension. 

 Some additional assumptions on the dynamics  help to clarify the topology of the carrying manifold. Let  $G_*(M^n)$ be a  class of Morse-Smale diffeomorphisms on a closed smooth connected $n$-dimensional manifold $M^n$ such that for any $f\in G_*(M^n)$ an  $(n-1)$-dimensional invariant manifold of  any  saddle periodic point of Morse index one or $(n-1)$ either do not intersect  invariant manifolds of any  other saddles or intersect only one-dimensional invariant manifolds.     
 If $W^u_p\cap W^s_q=\emptyset$ for any pair of saddle periodic points of a Morse-Smale diffeomorphism $f$,  we will say that $f$ {\it has no heteroclinic intersections}.  We set  
 \begin{equation}
     g_{f}=(k_1+k_{n-1}-k_0-k_n+2)/2. \end{equation}

 For $g\geq 0$  denote by $\mathbb S^n_g$ a connected sum of the $n$-dimensional sphere $S^n$ and  $g$ copies of the direct product $S^{n-1}\times S^1$. 

In~\cite{BGMP}, the following result is proved for $n=3$. 

\begin{stat}\label{orient-3} Let $f\in G_*(M^3)$ and $M^3$ be orientable. Then $g_{f}\geq 0$ and $M^3$  is diffeomorphic to $\mathbb S^3_{g_f}$. Moreover, for any $g \geq  0$ the  manifold $\mathbb S^3_g$  admits a diffeomorphism $f\in G_*(\mathbb S^3_g)$.
 \end{stat}

A series of papers~\cite{GrGuPo15, GrGu22_en, GrGu-Sib, GrGu23_en}, ~\cite{GrZhMe17},~\cite{GrMeZh22}  allows to generalize Statement~\ref{orient-3} for $n\geq 4$ as  follows. Recall that a Morse-Smale diffeomorphism is called {\it polar} if its non-wandering set consists of exactly one sink, one source and finite number of saddle periodic points. 

 \begin{stat}\label{orient-n+}  Let $M^n$ be orientable,  $n\geq 4$, and $f\in G_*(M^n)$.  Then 
 
 \begin{enumerate}
     \item 
     $g_{f}$ is a non-negative integer, and  $M^n$ is homeomorphic to a connected sum of  $\mathbb S^n_{g_{f}}$ and a  simply connected closed manifold $N^n$;
     \item $N^n$ admits a polar diffeomorphism $f\in G_*(N^n)$  without saddle periodic points of Morse indices $1$ and $(n-1)$; 

     \item if invariant manifolds of different saddle periodic points of $f$ of Morse indices  $i\in \{2,\dots, n-2\}$ do not intersect each other, then     $k_2=\dots=k_{n-2}=0$ if and only if $N^n$ is homeomorphic to the sphere $S^n$.
     \end{enumerate}
  \end{stat}

A problem of the realisation of the system  on given manifolds was  partially solved in~\cite{GrGu-Sib}, \cite{GrGu23_en}, \cite{GuSa24-mz_en}.

In~\cite{OsPo},  a generalization of Statement~\ref{orient-3} for non-orientable case has been obtained. We provide a generalization of Statement~\ref{orient-n+}. Moreover, we show that the  statement does not depend on a smooth structure of $M^n$. Hence, it   may be extended  to   a natural topological analogue of  the class of Morse-Smale systems, namely, the dynamical  systems  with finite and hyperbolic  chain recurrent set, that we call, following~\cite{MePoZi, PoZi_RCD},  {\it regular systems}. The chain recurrent set of a regular  homeomorphism consists of a  finite number of topologically hyperbolic  periodic points (see~Proposition~\ref{ch-recc-hyperbolicity}). Stable, unstable invariant manifolds and Morse index of  a topologically hyperbolic periodic point is defined similar to the invariant manifolds and the  Morse index  of the hyperbolic periodic point.   We denote by $G(M^n)$  a class of regular  homeomorphisms such that an   $(n-1)$-dimensional invariant manifold of an arbitrary periodic point $p$ of Morse index  $(n-1)$ either do not intersect any  invariant manifolds of other saddle periodic points or intersect only one-dimensional invariant manifolds.  If $M^n$ is smooth, then  $G_*(M^n)\subset G(M^n)$.  
 \begin{theo}\label{main+}
Let $M^n$ be a topological  closed manifold,  $n\geq 4$ and $f\in G(M^n)$. Then the following alternatives holds.  

\begin{enumerate}
\item If $M^n$ is orientable, it is homeomorphic to a connected sum of $\mathbb S^n_{g_f}$ and a simply connected  closed manifold $N^n$.

\item If $M^n$ is non-orientable, it is  homeomorphic to a connected sum of  $\mathcal S^n_{g_f}$ of $g_{f}>0$ copies of non-trivial $S^{n-1}$-bundle over $S^1$, and    a  simply connected  closed manifold $N^n$.

\item if invariant manifolds of different saddle periodic points of $f$ of Morse indices  $i\in \{2,\dots, n-2\}$ do not intersect each other, then     $k_2=\dots=k_{n-2}=0$ if and only if $N^n$ is homeomorphic to the sphere $S^n$.
     \end{enumerate} 
\end{theo} 

If  $f^t$ is a regular  flow  such that all $(n-1)$-dimensional invariant manifolds of its saddle equilibrium  states do not intersect any  other invariant manifolds of saddle equilibria,  then  a time-one shift map $f=f^1$ belongs to $G(M^n)$.  Hence,   all statements above are  also true for $f^t$.  

Before the proof of Theorem~\ref{main+} we provide  in Section~\ref{sec-con-sum}  an accurate proof of the fact that the connected sum of non-orientable topological manifolds is well-defined for dimension $n\geq 4$ that we could not find in a literature. In Theorem~\ref{trap} we prove that  regular homeomorphism, alike  their smooth prototypes, have a fine filtration, that may have an independent interest. 

\medskip 
{\bf Acknowledgements.} 
Authors thank to participants of the seminar of Laboratory of Dynamical systems and applications in NRU HSE for motivation and  useful discussions. 

\medskip 
{\bf Funding.} The work is supported by the Russian Science Foundation (project 23-71-30008),   except Sections~\ref{orient}-\ref{sec1.4},\ref{fin} that  were prepared within the framework of the project "International academic cooperation" of HSE University. 

\section{Definitions and auxiliary results}
\subsection{Notations}
\begin{itemize}
    \item $\mathbb R^n$  denotes the Euclidean space of dimension $n\geq 1$. For $k < n$ the space $\mathbb R^k$ is considered as a subset of $\mathbb R^n$ determined by condition $x_{k+1} = \dots = x_n = 0$; and  $\mathbb R^k_+$ is a subset of $\mathbb R^k$ determined by the inequality $x_k \geq 0$.

    \item $M^n$ is a closed topological manifold of dimension $n \geq 1$.

    \item $\partial M$ is a boundary of a  manifold $M$.
\item $\Int M$ is an interior of a manifold $M$.

\item {$M\cong N$ means that manifolds $M,\, N$ are homeomorphic}.

\item $S^{n-1}$, $B^{n}$, $n\geq 1$, denote  the topological  $(n-1)$-dimensional sphere and the  $n$-dimensional compact ball, that are  manifolds homeomorphic  to $$\mathbb S^{n-1}=\{x\in \mathbb R^{n}: ||x||=1\}, \mathbb B^n=\{x\in \mathbb R^n: ||x||\leq 1\},$$ correspondingly. $A^n,\, n \geq 1$, is an annulus of dimension $n$, that is a  topological manifold, homeomorphic to  $S^{n-1} \times [0,1]$. 

\item $\mathbb{CP}^2$ is  a complex projective plane, that is  a smooth four-dimensional manifold that is a  factor space  of $\mathbb{C}^3 \setminus \{O\}$ by the equivalence relation: $(z_1,z_2,z_3) \sim (\lambda z_1, \lambda z_2, \lambda z_3)$,  $\lambda \in \mathbb{C} \setminus \{0\}$. We suppose that $\mathbb {CP}^2$  is enriched  with an   orientation induced by the canonical  orientation of $\mathbb{C}^3$.

\item $\mathbb H^{n}_i= \mathbb B^i\times \mathbb B^{n-i}$ is an {\it $n$-dimensional $i$-handle}, $\mathbb F^{n-1}_i = \partial \mathbb B^i \times \mathbb B^{n-i}$ is a {\it  foot of $\mathbb H^{n}_i$}. 

\item $\Id_{_X}$ is an identity map on the set $X$.

\item $H_k(M)$ is  a $k$-dimensional singular homology group of the manifold $M$ with integral coefficients.

\item $H_k(M, B)$ is a  group of relative homology for a subspace  $B \subset M$  with integral coefficients.

\item $W^s_p, W^u_p$ are stable and unstable invariant manifolds of a hyperbolic periodic point $p$.

\item $W^s_{_X} (W^u_{_X})$ are unions of stable  (unstable) invariant manifolds of all hyperbolic periodic points from the set $X$.

\end{itemize}

\subsection{Exact sequences}\label{MV-pary}

Recall that a {\it direct sum} of two Abelian groups $G, H$ with binary operations $*, \times$, respectively, is the group $G \oplus H = \{(g, h): g \in G_1, h \in G_2\}$ with the binary operation $(g_1, h_1) \circ (g_2, h_2) = (g_1 * g_1, h_1 \times h_2)$ (see~\cite[Chapter 4.1]{Pr-hom_en}).

A   sequence of homomorphisms
\begin{equation}\label{exact-sequence}
\dots \to A_{n} \overset{\alpha_{n}}{\to} A_{n-1} \overset{\alpha_{n-1}}{\to} A_{n-2} \to \dots
\end{equation} is called {\it exact} if ${\rm Ker}\, \alpha_{n-1} = {\rm Im}\, \alpha_n$ for any $n$. 

\begin{prp}\label{exact-isomorph}
    Let the sequence~{\rm (\ref{exact-sequence})} be exact, $E_{n}: A_n \to B_{n}$ are  isomorphisms, and  $\beta_n = E_{n-1} \alpha_n E^{-1}_n$.  Then the sequence \begin{equation}
        \dots \to B_{n} \overset{\beta{n}}{\to} B_{n-1} \overset{\beta{n-1}}{\to} B_{n-2} \to \dots
    \end{equation} is exact.
\end{prp}
\begin{proofn}
    Since the sequence~(\ref{exact-sequence}) is exact and $E_n$ is an isomorphism, the following equality holds 
    \begin{equation}
    {\rm Im}\, \beta_n = E_{n-1} \alpha_n E^{-1}_{n}(B_n) = E_{n-1} \alpha_n(A_n) = E_{n-1}({\rm Im}\, \alpha_n) = E_{n-1}({\rm Ker}\, \alpha_{n-1}).
    \end{equation} Let $a \in {\rm Ker}\, \alpha_{n-1}$. Then  $a' = \alpha_{n-1}(a) = 0$ and $E_{n-2}(0) = 0$. So, $a \in {\rm Ker}\, E_{n-2} \alpha_n$. Since  $E_{n-1}$ is a one-to one correspondence,  $E_{n-1}({\rm Ker}\, \alpha_{n-1}) = {\rm Ker}\, \beta_{n-1}$. Then ${\rm Im}\, \beta_n = {\rm Ker}\, \beta_{n-1}$.
\end{proofn}
  
The following statement is a famous Mayer-Vietoris theorem (~\cite{Ma29}, \cite{Vi}, see also~\cite[Section 2.2]{Ha}).

 \begin{stat}\label{M-V}
 Let $X$ be a topological space and $U,\, V$ be subspaces of $X$ such that $X = {\rm int}\, U \cup {\rm int}\, V$. Set $N = U \cap V$. Then the  sequence \begin{equation}\label{mvseq}
 \dots \to H_k(N) \overset{i_*}{\to} H_k(U) \oplus H_k(V) \overset{j_*}{\to} H_k(X) \overset{\partial_*}{\to} H_{k-1}(N) \to \dots,
 \end{equation} where $i_*([z]_{_N}) = ([z]_{_U}, [-z]_{_V})$; $j_*([x]_{_U}, [z]_{_V}) = [x+z]_{_X}$; $\partial_*([z]_{_X}) = [\partial\, z_{_U}]_{_T} = [-\partial\, z_{_V}]_{_N}$, where $z_{_U} = z \cap U$, $z_{_V} = z \cap V$, is exact.
  \end{stat}

  The sequence determined in  Statement~\ref{M-V} is called a {\it Mayer-Vietoris sequence}. Another classic  instrument in algebraic topology that we use below,   is an exact sequence of a pair of topological spaces (see~\cite[Theorem 2.16]{Ha}).

  \begin{stat}\label{paryseq}
Let $X$ be a topological space and $A$ be its subspace. Then there is an exact sequence 
\begin{equation}
\dots \to H_k(A) \overset{i_*}{\to} H_k(X) \overset{j_*}{\to} H_k(X,A) \overset{\partial_*}{\to} H_{k-1}(A) \to \dots,
\end{equation} where $i_*([a]_{_A}) = [a]_{_X}$, $j_*([x]_{_X}) = [x]_{_{(X,A)}}$, $\partial_*([x]_{_{(X,A)}}) = [\partial x]_{_{A}}$.
  \end{stat}

\subsection{\normalsize{Embedding in topological manifolds and topological transversality}}

A Hausdorff second-countable space $M$ is called a {\it topological $n$-dimensional manifold} (or a {\it topological $n$-manifold}) if every point $x \in M$ has a neighbourhood homeomorphic to $\mathbb R^n$ or to $\mathbb R^n_+$. {\it A  boundary of $M$} is a  set $\partial M$ of points in $M$ that  have  neighbourhoods homeomorphic to $\mathbb R^n_+$. {\it An interior of $M$} is a set  $\Int M$  of points that have  neighbourhoods homeomorphic to $\mathbb R^n$. 

A  map  $e: X\to Y$  of topological space is called {\it the topological embedding} if $e: X\to e(X)$ is a homeomorphism, where $e(X)$ is considered with the topology induced by the topology of $Y$.  

{Recall that an {\it isotopy} of a topological space  $Y$ is a continuous map  $H: Y\times [0,1] \to Y$ such that for any $t \in [0,1]$ a  map $h_t(x)=H(x, t)$ is a homeomorphism. If $X \subset Y$ and $h_t|_{X} = \Id_{_X}$ for any $t \in [0,1]$ then the  isotopy is called  {\it  relative to $X$}.} 

Topological embeddings $e_1, e_2: X\to Y$ are called {\it ambient isotopic} if there exists an isotopy $h_t: Y\to Y$ such that $h_0=\Id_{_Y}$
and $e_2=h_1e_1$.

A manifold  $X \subset M$ with possibly a non-empty boundary is called   {\it  a submanifold of $M$ or locally flat  in $M$} if for any point  $x \in \Int X\  (x\in \partial X)$ there is a neighbourhood $U_x\subset M$  and  a homeomorphism $\psi: U_x\to \mathbb R^n$ such that $\psi(U_x \cap X) = \mathbb R^k$ ($\psi(U_x \cap X) = \mathbb R^k_+$). 

A submanifold $X\subset M$ is {\it proper} if either $\partial X=\emptyset$ and $X\in \Int M$ or   $\partial X\subset \partial M$ and $\Int X\subset \Int M$.
 
Let $X\subset  M^n$ be a closed  manifold of dimension $(n-1)$. If   there exists an embedding $e: X \times [0,1) \to M$ such that $e(X \times \{0\}) = X$, then the set  $C = e(X \times [0, 1))$ is called {\it a collar} of $X$. A boundary $\partial M^n$ always have a collar. If there is an embedding $e: X \times [-1,1) \to M$ such that $e(X \times \{0\}) = X$ then the manifold $X$  is called bi-collared. According to~\cite{Br}, a boundary $\partial M^n$ of a compact topological manifold have a collar and any two-sided locally flat $(n-1)$-dimension manifold in $M^n$ is bi-collared. It is well-known that if $M^n$ is  orientable, then any orientable locally flat   $(n-1)$-dimensional manifold $X\subset M^n$ is two-sided. Hence,    $X$ is bi-collared.  In~\cite[Proposition 2]{GuSa24-pde} the following statement is obtained. 

\begin{stat}\label{bicol} Let $M^n$  be a closed topological manifold (either orientable or not) and  $S^{n-1} \subset M^n$ be a locally flat sphere. Then $S^{n-1}$  is bi-collared. 
    \end{stat}

In other words,  Statement~\ref{bicol} means that a locally flat sphere $S^{n-1}$ have a topological analog of  tubular neighbourhood in $M^n$. For arbitrary topological manifolds, a generalisation of the notions of   tangent bundle and  the  tubular neighbourhood involves a following notion of  microbundle introduced by J.~Milnor in~\cite{M}. 

A $k$-microbundle over a topological space  $B$ is a set $\{E, B, p, i\}$, where $E$ is a topological space and $p: E\to B, i: B\to E$ are continuous maps such that $pi|_{_B}=\Id_{_B}$ 
and for every $b\in B$  there are  neighbourhoods   $U\subset B, V\subset E$ of $b, i(b)$, respectively,  and  a homeomorphism $\gamma:V\to U\times \mathbb R^k$ such that such that $i(U)\subset V, p(V)\subset U$, and   $\gamma i|_{_U}=i'$, $p'\gamma|_{_{V}}=p|_{_V}$, where $i'(x)=(x,0)$ is an inclusion of $U$ to $U\times \mathbb R^k$ and $p' (x,y)=(x,0)$ is a projection of $U\times \mathbb R^k\to U\times \{O\}$.

Spaces $E, B$ are  called {\it  the total space and the base}  and $i, p$ are called the  {\it injection and projection} of the microbundle, respectively.   The set $\{B\times \mathbb R^k, B,  p', i'\}$ is a natural example of the microbundle that is called {\it the trivial microbundle}. 

Let $M$  be a topological manifold,  $X\subset M$  be a submanifold, $i_{_X}: X\to M$ be an inclusion. $X$ is said to have {\it a normal $k$-microbundle in $M$} if   there exists a neighbourhood $E\subset M$ of $X$ and a retraction $p: E\to X$ such that $\xi^k=\{E, X, p, i_{_X}\}$ is a microbundle. 

A manifold $Y\subset M$ is called {\it transversal to a normal microbundle $\xi^k$ in $M^n$ or embedded  topologically transversely to $X$} if $Y\cap X$  is a submanifold of $Y$ with a normal microbundle $\nu^k$ in $Y$ which is a restriction of $\xi^k$, that is an inclusion $i_{_{Y}}: Y\to M$ induces an open topological embedding of each fiber of $\nu^k$  to a some fiber of $\xi^k$.  Let us remark that if $Y$ is  embedded  topologically transversely to $X$, then $k=\Dim Y-\Dim(X\cap Y)=\Dim M-\Dim X$, so  
\begin{equation}
\Dim (X\cap Y)=\Dim X + \Dim Y - \Dim M.
\end{equation}

$Y$ is embedded topologically  transversely to $X$  {\it near $C\subset X$} if there is an open neighborhood $M_0\subset M^n$ of $C$ such that $M_0\cap Y$ is embedded topologically transversely to $M_0\cap X$ in $M_0$. 

According to~\cite[Theorem~1.5 of Essay III]{KSi} and~\cite{Quinn}, the following statement is true.

\begin{stat}\label{quinn} Let $X, Y\subset M$ be  proper  compact submanifolds,  $X$ has a normal microbundle in $M$ and there are closed subsets $C\subset D\subset M$ such that $Y$ is embedded topologically transversal to $X$  near $C$. Then there is an  isotopy $h_t: M\to M$  supported in any given neighbourhood of $(D \setminus C)\cap X\cap Y$ such that $h_0=id|_{_M}$ and $h_1(Y)$ is embedded topologically transversal to $X$ near $D$.  
\end{stat}

\subsection{Orientations of topological manifolds}\label{orient}

We recall   notions of a local and a global orientation of a topological $n$-dimensional connected manifold $M$ without boundary, $n \geq 2$, that do not depend on the  existence  of PL or smooth structures on $M$. Mainly we follow~\cite[Section 3.3]{Ha} but clarify some  details that will be  used in Section~\ref{sec-con-sum} in a  proof of the fact  that a connected sum of  topological  non-orientable  manifolds is well-defined.

Let $B \subset M$ be  a compact $n$-dimensional ball locally flat in $M$. It follows from~\cite[Theorem 2]{Br} that  the sphere $\partial B = S$ is {\it collared in $M \setminus \Int B$}, that is there exists an  embedding $e: S \times [0,1) \to M \setminus \Int B$ such that $e(S \times \{0\}) = S$. Set $C = e(S \times [0, 1))$. Let $x \in \Int B$. The set $U_x = B \cup C$ is a neighbourhood of $x$ homeomorphic to $\mathbb R^n$ and the inclusion $(U_x, C)\to (M, M \setminus {\rm int}\, B) $ induces a  homomorphism $i_{_B}: H_n(U_x, C) \to H_n(M, M \setminus {\rm int}\, B)$. Due to  Excision Theorem (\cite[Theorem 2.20]{Ha}),  $i_{_B}$   is an isomorphism. According to Statement~\ref{paryseq} there is the exact sequence \begin{equation}
    \dots \to H_n(U_x) \to H_n(U_x, C) \overset{{\partial_{_B}}}{\to} H_{n-1}(C) \to H_{n-1}(U_x) \to \dots
\end{equation} Since the groups $H_{n}(U_x),\, H_{n-1}(U_x)$ are trivial for $n \geq 2$, the map $\partial_{_B}: H_n(U_x, C) \to H_{n-1}(C)$ is an isomorphism. The manifolds $C$ and $S$ are homotopy equivalent, then by~\cite[Corollary 2.11]{Ha} there is an isomorphism $F_*$ between $H_{n-1}(C)$ and  $H_{n-1}(S)$.  Define an isomorphism  \begin{equation}
I_*: H_n(M, M \setminus {\rm int}\, B) \to H_{n-1}(S)
\end{equation} setting $I_*=F_*\partial_{_B}i_{_B}^{-1}$. 
The group $H_{n-1}(S)$ is infinite cyclic (see, for instance,~\cite[Example 2.5]{Ha}).  A generator $\mu_{_S}$ for $H_{n-1}(S)$ is called a {\it fundamental class of $S$}. Below we show that $\mu_{_S}$ determines an orientation of $S$. The  group $H_n(M, M \setminus {\rm int}\, B)$ is infinite cyclic, too, and $I_*^{-1} (\mu_{_S}) = \mu_{_B}$ is its generator.

Since an inclusion  $M\setminus \Int B\to M\setminus \{x\}$ is a homotopy equivalence, the map 
\begin{equation}
j_x: H_n(M, M \setminus {\rm int}\, B) \to H_n(M, M \setminus \{x\})
\end{equation} induced by the inclusion, is an  isomorphism. We  call  $j_x$ a {\it natural isomorphism}. As a corollary of all facts above  we immediately get the following proposition.

\begin{prp}\label{loc-or}
For any $x \in M$ the group $H_n(M, M \setminus \{x\})$ is infinite cyclic.
\end{prp}

The generator $\mu_x$ of $H_n(M, M \setminus \{x\})$ is called a {\it local orientation of $M$ at the point $x$}. There are exactly two local orientations of $M$ at any point $x$: $\mu_x$ and its opposite $-\mu_x$.

Let $\mu_{x'}$ be the local orientation at a point  $x'\in {\rm int}\, B \setminus \{x\}$. Local orientations $\mu_x,\, \mu_{x'}$ are  {\it locally consistent} if $j_{x}^{-1}(\mu_x)=j_{x'}^{-1}(\mu_{x'})=\mu_{_B}$. The manifold $M$ is called {\it orientable} if {it is possible to choose local orientation $\mu_x$ at every point $x \in M$} such that for any locally flat compact ball $B\subset M$ and any  $x,\, x'\in {\rm~int}~B$,  local orientations $\mu_x,\, \mu_{x'}$ are  locally consistent.

Let  $\widetilde M$ be a set of all possible local orientations of points $x\in M$.  
The set  $\widetilde M$ can be enriched by a topology  $\widetilde \tau$ with a basis formed by the sets   $U(\mu_{_B}) = \{\mu_x: x \in {\rm int}\, B, \mu_x = j_x (\mu_{_B})\}$, where  $B\subset M$ is  a compact locally flat $n$-ball,   $\mu_{_B}$  be a generator of  $H_n(M, M \setminus {\rm int}\, B)$ and $j_x: H_n(M, M \setminus {\rm int}\, B) \to H_n(M, M \setminus \{x\})$ is the natural isomorphism.
Then a map   $p: \widetilde M \to M$   that  puts in correspondence to each point  $\mu \in \widetilde M$  a point $x\in M$ such that $\mu_x=\mu$,   is a two-fold covering map. 

Due to~\cite[Proposition 3.25]{Ha}, the following statement holds. 

\begin{stat}\label{connected-cover} $\widetilde M$  is a topological $n$-dimensional orientable manifold. $\widetilde M$ is connected if and only if $M$ is non-orientable. 
\end{stat}

According to~\cite[Theorem 3.26]{Ha}, the following statement holds.

\begin{stat}\label{Ha3}
Suppose $M$ to be closed. If $M$ is orientable then for any $x \in M$,  the inclusion map of $(M, \emptyset)\to (M, M \setminus \{x\})$ induces an isomorphism $\Delta_x: H_n(M) \to H_n(M, M \setminus \{x\})$. Hence, $H_n(M)$ is infinite cyclic.  If $M$ is non-orientable, then $H_n(M)$ is trivial.  
\end{stat}

\begin{rem}\label{j_{_B}}
    Let $j_{_B}: H_n(M) \to H_n(M, M \setminus {\rm int}\, B)$ be a homomorphism induced by the inclusion $(M, \emptyset) \to (M, M \setminus {\rm int}\, B)$.  Then the  isomorphism $\Delta_x$ defined in Statement~\ref{Ha3} is a composition of $j_{_B}$ and the  natural isomorphism $j_{x}: H_{n}(M, M \setminus {\rm int}\, B) \to H_n(M, M \setminus \{x\})$.  Hence,  $j_{_B}$ is an isomorphism.
\end{rem}

For the  closed orientable manifold $M^n$, the generator  $\mu_{_{M^n}}$ is called  a {\it fundamental class}, and a  fixed fundamental class $\mu_{_{M^n}}$ is called an {\it  orientation of $M^n$}. The image $\Delta_x(\mu_{_{M^n}})$ can be considered as the  local orientation of $M^n$ at a point $x$. Due to Remark~\ref{j_{_B}}, for any ball $B$ and any two points $y,\, y' \in {\rm int}\, B$ the local orientations $\Delta_y(\mu_{_{M^n}}),\, \Delta_{y'}(\mu_{_{M^n}})$ are locally consistent. On the other hand, locally consistent orientations of points of $M^n$ also determine the same generator $\mu_{_{M^n}}$. Moreover, for any compact locally flat $n$-ball $B\subset M^n$,   the isomorphism $I_* j_{_B}: H_n(M^n) \to H_{n-1}(S)$ induces the orientation $\mu_{_S} = I_* j_{_B} (\mu_{_{M^n}})$ of the sphere $S$ that will be called a {\it natural orientation of $S$ induced by $\mu_{_{M^n}}$}.

\begin{prp}\label{minus-ball} Let $M^n$ be non-orientable manifold and   $X\subset M^n$ be a compact locally flat $n$-ball. Then a manifold   $M' = M^n \setminus X$ is non-orientable. \end{prp}
{\begin{proofn} Set $S=\partial X$. Since $S$ is locally flat sphere, due to Statement~\ref{bicol} it is bi-collared, that is there exists a topological embedding  $e:\mathbb S^{n-1}\times [-1,1]\to M^n$ such that  $S=e(\mathbb S^{n-1}\times \{0\}$. Set   $C_{_X}=e(\mathbb S^{n-1}\times [-1,1])\cap X$,  $C' =e(\mathbb S^{n-1}\times [-1,1])\cap M'$,  $Y = X \cup C'$, $M'' = M' \cup C_{_X}$. According to Statement~\ref{M-V},  the following sequence is exact:  \begin{equation}\label{MV}
    0 \to H_n(M^n) \overset{\partial_*}{\to} H_{n-1}(Y \cap M'') \overset{i_*}{\to} H_{n-1}(Y) \oplus H_{n-1}(M'') \overset{j_*}{\to} H_{n-1}(M^n) \to \dots\end{equation} The set $Y \cap M''$ is homotopy equivalent to $S$, $M''$ is homotopy equivalent to $\Cl M'$,  $Y$ is homotopy equivalent to $X$, and $H_n(M^n) = 0$ since $M^n$ is non-orientable. By Proposition~\ref{exact-isomorph} the sequence~(\ref{MV}) transforms to there is the exact sequence \begin{equation}\label{MV1} H_n(M^n) = 0 \overset{\Delta*}{\to} H_{n-1}(S) \overset{I_*}{\to} H_{n-1}(X) \oplus H_{n-1}(\Cl M') \overset{J_*}{\to} H_{n-1}(M^n) \to \dots \end{equation} $H_{n-1}(S)$ is isomorphic to $\mathbb Z$. Suppose $M'$ is orientable. Then ${\rm cl} M'$ is an $(n-1)$-chain bounded by the cycle $S$. Hence, the class $[S]$ is  trivial in $H_{n-1}({\rm cl} M')$,  ${\rm Im} I_* = 0$  and the sequence above means that the sequence 
    \begin{equation}
    0 \to \mathbb Z \to 0
    \end{equation} is exact, {that is impossible}. Then $M'$ is non-orientable. \end{proofn}}

Let $N^n$ be an oriented  manifold with an orientation $\mu_{_{N^n}}$,   $h: M^n \to N^n$ be a homeomorphism,  and $h_*: H_n(M^n) \to H_n(N^n)$  be an isomorphism induced by $h$. If $h_*(\mu_{_{M^n}}) = \mu_{_{N^n}}$, then $h$ is called {\it orientation  preserving}, otherwise   $h$  is called {\it  orientation reversing}. 

\subsection{Connected sums of closed topological manifolds}\label{sec-con-sum}

Let $X,\, Y$ be compact manifolds,  $A \subset \partial X,\, B \subset \partial Y$ be submanifolds,  and  $\xi: A \to B$ be a homeomorphism. Then   a factor space $X\cup_\xi Y$ of   $X \sqcup Y$ by a minimal equivalence relation $\sim$ such that $a \sim \xi(a)$,  is a manifold  said to be {\it obtained by gluing $X$ to $Y$ by means of $\xi$}. For a point $p \in X \sqcup Y$ we denote by $[p]_{\xi}$ the equivalence class of $p$ with respect to  this equivalence relation.

Let $B^n_{_M},\, B^n_{_N}$ be locally flat balls in closed manifolds $M^n,\, N^n$, respectively, $S^{n-1}_{_M} = \partial B^n_{_M},\, S^{n-1}_{_N} = \partial B^n_{_N}$, and  $\xi: S^{n-1}_{_M} \to S^{n-1}_{_N}$ be  a homeomorphism. 
 If one of $M^n,\, N^n$ is  non-orientable, then we call  a manifold $M^n \#_{\xi} N^n$  obtained by gluing $M^n \setminus {\rm int}\, B^n_{_M}$ and $N^n \setminus {\rm int}\, B^n_{_N}$ by means  of  $\xi$, {\it a  connected sum} of $M^n,\, N^n$.  
 If $M^n$ and $N^n$ are orientable, then  we fix their orientations $\mu_{_M},\, \mu_{_N}$ and natural orientations $\mu_{S_{_M}},\, \mu_{S_{_N}}$  of the spheres $S^{n-1}_{_M},\, S^{n-1}_{_N}$, and   assume that $\xi$  reverses the orientations. Then orientations $\mu_{_M},\, -\mu_{_N}$ determine an orientation on a   manifold $M^n \#_{\xi} N^n$  obtained by gluing $M^n \setminus {\rm int}\, B^n_{_M}$.  We call this manifold {\it a connected sum of $M^n,\, N^n$}.   Let   $\xi_*: H_{n-1}(S_{1}) \to H_{n-1}(S_{2})$ is   an  isomorphism induced by $\xi$. 
     
 If  $M^n,\, N^n$ are   smooth  (or PL)  oriented manifolds, and $\xi$ is an orientation-reversing diffeomorphism (PL-homeomorphism),  then, according to~\cite[Lemma 2.1]{KerMil} (\cite[Chapter 3]{RuSa_en}),  the connected sum of  $M^n,\, N^n$ is well-defined, that is  does not depend on a choice of $B^n_{_M},\, B^n_{_N}$ and  $\xi$. In{~\cite{BrGl}} it is shown that the connected sum, without any restriction on the gluing homeomorphism,  is  well-defined  for  topological manifolds if at least one of the summand is  {\it homogeneous}. According to~\cite{BrGl}, a manifold $M^n$  is called homogeneous if  for any locally flat embeddings $e_1,e_2:B^n\to M^n$ there exists a homeomorphism $h:M^n\to M^n$ such that $e_2=he_1$. In this section we show  that the connected sum is well-defined for arbitrary closed topological manifolds $M^n,\, N^n$  that may be non-smoothable,  non-triangulable, and   non-orientable.  In fact,  in~Corollary~\ref{loop2} and Proposition~\ref{hempel}  we show  that all non-orientable manifolds as well as oriented manifolds that admit orientation reversing homeomorphisms  are homogeneous. 

The following theorem is the summary of this section.
 \begin{theo}\label{group-operation}
   The connected sum of topological  manifolds    $M^n,\, N^n$ does not depend on the choice of the  balls $B^n_{_M},\, B^n_{_N}$ and the gluing map  $\xi$, and the  following properties hold:
   \begin{enumerate}
        \item $M^n \# N^n \cong N^n \# M^n$;
        \item $M^n \# S^n \cong M^n$;
        \item $(M^n \# N^n) \# L^n \cong M^n \# (N^n \# L^n)$.
   \end{enumerate}   
\end{theo}

\begin{rem}\label{lem2-rema}
Due to  Theorem~\ref{group-operation} we  may omit the mention of the  gluing map in the definition of the connected sum and denote it by $M^n \# N^n$. We will denote a  connected sum of $m \geq 0$ copies of the manifold $M^n$  by $m M^n$, setting $0M^n = S^n$. 
 \end{rem}

The proof of  Theorem~\ref{group-operation} is given in Propositions~\ref{balls-iso}-\ref{consum-nonor} below. 
  
\begin{prp}\label{balls-iso}
    For any compact $n$-dimensional submanifold $X$ of the $n$-manifold $M$ and for any two compact $n$-balls $B,\, B' \subset M \setminus X$ locally flat in $M$ there exists an  isotopy $h_t: M \to M$ relative to $X$ such that $H_0 = {\rm id}_{_M}$ and $h_1(B)=B'$. 
\end{prp}
 
\begin{proofn} We prove the proposition in three steps.

{\bf Step 1.} Suppose that $B'\subset \Int B$ and construct the desired isotopy. The spheres $S' = \partial B',\, S = \partial B$ are locally flat in $M$. It follows from the Annulus Theorem (see~\cite[Section 14.2]{GrGu-book_en} for references) that the domain in $M$ bounded by $S'$ and $S$ is homeomorphic to the annulus $S^{n-1} \times [0,1]$. Moreover, the spheres $S',\, S$ have collars $C',\, C$ in $B'$ and $M \setminus (\Int B \cup X)$, respectively. Let $\Sigma' = \partial C' \setminus S'$ and $\Sigma = \partial C \setminus S$. According to~\cite[Proposition 14.2]{GrGu-book_en} there is an embedding $e: S^{n-1} \times [0,1] \to M$ such that $e(S^{n-1} \times \{0\}) = \Sigma'$, $e(S^{n-1} \times \{s'_0\}) = S'$, $e(S^{n-1} \times \{s_0\}) = S$ and $e(S^{n-1} \times \{1\}) = \Sigma$, where $0 < s'_0 < s_0 < 1$. Set $K = e(S^{n-1} \times [0,1])$. Let $L: [0,1] \to [0,1]$ be a linear function determined by the formula \begin{equation}\label{otrezok}
    L(s) = \begin{cases}
\dfrac{s'_0}{s_0} s,\, s \in [0, s_0];\\
\dfrac{1-s'_0}{1 - s_0} (s - s_0) + s'_0,\, s \in [s_0, 1].
    \end{cases}
\end{equation} By definition $L$ is a self-homeomorphism of the segment $[0,1]$ such that $L(0) = 0$, $L(s_0) = s'_0$ and $L(1) = 1$. For  $t \in [0,1]$ we define a function $L_t: [0,1] \to [0,1]$ by the formula \begin{equation}\label{otrezok-iso}
    L_t(s) = t L(s) + (1-t) s.  
\end{equation} For every $t$ the function $L_t$ is a homeomorphism and $L_0 = \Id_{[0,1]}$, $L_t(0) = 0$, $L_t(1) = 1$,  $L_1 = L$. Let us define an isotopy $\phi_t: K \to K$ by the formula $\phi_t(e(p, s)) = e (p, L_t(s))$, where $p \in S^{n-1},\, s \in [0,1]$. By definition $\phi_0 = \Id_{_K}$, $\phi_1|_{_{\partial K}} = \Id_{_{\partial K}}$ and $\phi_1(S) = S'$. Hence, the isotopy $\phi_t$ extends to the isotopy $\Phi_t$ determined by the formula \begin{equation}\label{Phi}
    \Phi_t(x) = \begin{cases}
        \phi_t(x),\, x \in K;\\
        x,\, x \in M \setminus \Int K.
    \end{cases}
\end{equation} By construction, $\Phi_t$ is the required isotopy. 

{\bf Step~2.} Suppose that ${\rm int}\, B'\cap {\rm int}\,B \neq \emptyset$ and construct the desired isotopy. By condition, there exists a compact locally flat $n$-dimensional ball $B_0 \subset {\rm int}\, B\cap {\rm int}\, B'$. According to Step 1 there exists isotopies  $\Phi_t,\, \Phi'_t: M\to M$ relative to $X$ such that $\Phi_0=\Phi'_0={\rm id}_{_M},\, \Phi_1(B)=B_0,\, \Phi'_1(B')=B_0$. Then the map 
\begin{equation}\label{st2}
G_t(x)=\begin{cases}
   \Phi_{2t}(x), t\in [0,1/2];\\
   \Phi'^{-1}_{2t-1}(\Phi_1(x)), t\in [1/2,1]
\end{cases} 
\end{equation}
is the desired isotopy.  

{\bf Step~3.} Suppose that ${\rm int}\, B'\cap {\rm int}\,B=\emptyset$ and construct the desired isotopy. 
 Let $x \in {\rm int} B,\, x' \in {\rm int} B'$. Then by the Homogeneity Theorem (see~\cite[Lemm~3.33 of Chapter 3]{RuSa_en}) there is an isotopy $\widehat F_t: M \setminus \Int X \to M \setminus \Int$X relative to $\partial X$ such that $\widehat F_0 = {\rm id}_{_M}$ and $\widehat F_1(x) = x'$. The isotopy $\widehat F_t$ naturally extends to the isotopy $F_t: M \to M$ relative to $X$. Then $F_1({\rm int}\,B)\cap {\rm int}\, B'\neq \emptyset$ and due to Step~2 there exists an isotopy  $G_{t}: M\to M$ relative to $X$ such that $G_{0}={\rm id}_{_M}$ and $G_1(F_1(B))=B'$. Then the desired isotopy $h_t: M \to M$ is determined by the formula
\begin{equation}\label{st3}
h_t(x)=\begin{cases}
   F_{2t}(x), t\in [0,1/2];\\
   G_{2t-1}(F_1(x)), t\in [1/2,1].
\end{cases} 
\end{equation}
 \end{proofn}

Let $\widetilde B^{n}_{_M} \not= B^{n}_{_M}$ and $\widetilde B^{n}_{_N} \not=B^{n}_{_N}$  be locally flat $n$-balls in  $M^{n},\, N^{n}$, respectively, 
 and $\xi: \partial B^{n}_{_M}\to \partial B^{n}_{_N}$ be a homeomorphism. Due to   Proposition~\ref{balls-iso}  there exists homeomorphisms $h_{_M}: M^{n} \to M^{n}$,  $h_{_N}: N^{n} \to N^{n}$  such that $h_{_M}(B^{n}_{_M}) = \widetilde B^{n}_{_M}$ and $h_{_N}(B^{n}_{_N}) = \widetilde B^{n}_{_N}$. Define a homeomorphism $\tilde \xi: \partial \widetilde B^{n}_{_M} \to \partial \widetilde B^{n}_{_N}$
by $\tilde \xi = h_{_N} \xi h^{-1}_{_M}$. 

\begin{cor}\label{consum-balls}
     $M^{n} \#_{\xi} N^{n}$ is homeomorphic to $M^{n} \#_{\tilde \xi} N^{n}$.
\end{cor}
\begin{proofn}
Let $g: (M^{n} \setminus {\rm int}\, B^{n}_{_M}) \sqcup (N^{n} \setminus {\rm int}\, B^{n}_{_N}) \to (M^{n} \setminus {\rm int}\, \widetilde B^{n}_{_M}) \sqcup (N^{n} \setminus {\rm int}\, \widetilde B^{n}_{_N})$ be a homeomorphism determined by  \begin{equation}\label{con-sum-homeo}
        g(x) = \begin{cases}
        h_{_M}(x),\, x \in M^{n} \setminus {\rm int}\, B^{n}_{_M};\\
        h_{_N}(x),\, x \in N^{n} \setminus {\rm int}\, B^{n}_{_N},
        \end{cases}
    \end{equation}
    
    By definition of $\tilde \xi$ the following equalities hold for any $x\in M^{n} \setminus {\rm int}\, B^{n}_{_M}$:  
    \begin{equation}\label{xi'}
   \tilde \xi g(x) = \tilde \xi h_{_M}(x) = h_{_N} \xi h^{-1}_{_M} h_{_M}(x) = h_{_N} \xi(x) = g \xi(x).     
    \end{equation}

        Hence, $g([x]_\xi)=[g(x)]_{\tilde{\xi}}$. Similar property holds for any $x \in N^{n} \setminus {\rm int}\, B^{n}_{_N}$. Hence, $g$ induces a homeomorphism between $M \#_{\xi} N$ and $M \#_{\tilde \xi} N$.
\end{proofn}

\begin{prp}\label{loop}
     Let $M^n$ be a non-orientable closed manifold, $B^n, X \subset M^n$ be  locally flat balls,   $S^{n-1} = \partial B^n$. Then there is an isotopy $F_t: M^n \to M^n$ relative to $X$ such that:
     \begin{enumerate}
         \item $F_0 = {\rm id}_{_{M^n}}$;
         \item $F_1 (B^n) = B^n$; 
         \item $F_1|_{S^{n-1}}$ reverses an orientation of $S^{n-1}$.
     \end{enumerate}
 \end{prp}
\begin{proofn}
 Set $M'=M\setminus X$. {Due to Proposition~\ref{minus-ball}}, $M'$ is non-orientable. In Section~\ref{orient}, a two-fold covering map  $p:\widetilde M' \to M'$ is defined, where the covering space $\widetilde M'$ is the union of all possible local orientations of points $x\in M'$. By definition the point $x \in \Int B^n$ has a preimage $p^{-1}(x)$ consisting of two points $\mu_x,\, - \mu_x$. Due to  Statement~\ref{connected-cover},  $\widetilde M'$ is connected and then it is path connected. Hence, there is a path $\tilde \alpha: [0,1] \to \widetilde M'$ connecting $\tilde \alpha (0) = \mu_x$ with $\tilde \alpha (1) = -\mu_x$.  Set $\alpha=p \tilde \alpha$. Then $\alpha: [0,1] \to M'$ is a loop in $M'$ such that $\alpha(0)=\alpha(1)=x$. 

 Since the loop is compact, there exists a finite set  $0=t_0<t_1<t_2< \dots<t_k=1$ such that the segment  $\alpha([t_i, t_{i+1}])$ belongs to an open ball $B^n_{i+1} \subset M'$ properly covered by $p$. That means that there is a connected union of  balls $\{\widetilde B^n_{i+1}\}$  such that $\tilde \alpha([t_i,t_{i+1}]) \subset \widetilde B^n_i$ and $p|_{\widetilde{B}^n_i}: \widetilde{B}^n_i\to B^n_i$ is a homeomorphism.  Set $x_i=\alpha(t_i)$, $i\in \{0,\dots, k\}$.

 Due to Proposition~\ref{balls-iso},  without loss of generality we may assume that $B^n \subset B^n_1\cap B^n_{k}$.  Moreover, for any $i\in \{1,\dots, k\}$ there exists an isotopy $h^i_{t}:M'\to M'$ relative to $M'\setminus B^n_i$   such that  $h^i_0=\Id_{_{M^n}}$, $h^i_1(x_{i-1})=x_{i}$, and  $h^{k}_1 \cdots h^1_1 (B^n) = B^n$. Set $\beta_i=\bigcup\limits_{t\in [0,1]}h^i_t(x_{i-1})$, $\beta=\bigcup\limits_{i=1}^k\beta_i$. By definition, $\beta_i$ is a path connecting points $x_{i-1}, x_i$, and $\beta$ is a loop. By  uniqueness of the lifting (see, for instance,~\cite[Lemma 17.4]{Ko_en}), there is a unique path $\tilde\beta_i \subset \widetilde B^n_i$  such that $p(\tilde\beta_i)=\beta_i$ that begins at $\tilde\alpha_{t_{i-1}}$. The path $\tilde \beta = \bigcup\limits_{i=1}^k \tilde \beta_i$ is a lifting of $\beta$ that begins at $\mu_x$ and ends at  $-\mu_x$. Then the map $F'_t$ given by the formula: \begin{equation*}
     F'_t(z) = h^i_{kt-(i-1)}(z),\, t \in \left[\dfrac{i-1}{k}, \dfrac{i}{k}\right].
     \end{equation*} is the isotopy relative to $M' \setminus \bigcup\limits_{i = 1}^{k} B^n_i$ reversing the orientation of $S^{n-1}$. The isotopy $F'_t$ naturally extends to the required isotopy $F_t: M^n \to M^n$ relative to $X$.  
   \end{proofn} 
 \begin{cor}\label{loop2} Any non-orientable closed manifold $M^n$ is homogeneous.
 \end{cor}
 \begin{proofn} Let $e_1, e_2:\mathbb B^n\to M^n$ be locally flat embedding. Due to Proposition~\ref{balls-iso} we may assume that $e_1(\mathbb B^n)=e_2(\mathbb B^n)$.  Set $B^n=e_i(\mathbb B^n)$, $S^{n-1}=\partial B^n$ and denote by $C$ a collar of $S^{n-1}$ in $M^n\setminus \Int\, B^n$. We determine  any point $p\in C$ by two coordinates $x\in S^{n-1}, t\in [0,1]$.
 Two cases are possible: a map  $e_2e^{-1}_1|_{S^{n-1}}:S^{n-1}\to S^{n-1}$  preserves or  reverses an  orientation of $S^{n-1}$. 
 In the first case there exists an isotopy $g_t: S^{n-1}\to S^{n-1}$ such that $g_0={0}=e_2e^{-1}_1, g_1=\Id_{_{S^{n-1}}}$. Set 
 \begin{equation}
 h(p)=\begin{cases}
 e_2e^{-1}_1(p), p\in B^n;\cr
  (g_t(x), t), p=(x,t)\in C;\cr
  x, x\in M^n\setminus (B^n\cup C).
     \end{cases}
\end{equation}

 Then $he_1=e_2$, hence $M^n$ is homogeneous. 

 In case 2, due to Proposition~\ref{loop}, there exists an isotopy $F_t:M^n\to M^n$  such that $F_1(B^n)=B^n$ and  $F_1|_{S^{n-1}}$ reverses the orientation of $S^{n-1}$. Then $F_1e_2e_{1}^{-1}|_{S^{n-1}}$ is an  orientation preserving homeomorphism,   and  as above one may construct a homeomorphism $\eta: M^n\to M^n$  that coincides with $F_1e_2e_{1}^{-1}$ on $B^n$.   Then $h=F_{1}^{-1}\eta$  satisfy the condition $he_1=e_2$ and, hence, $M^n$ is homogeneous. 
 \end{proofn}
 
\begin{prp}\label{consum-nonor}
   If $M^n,\, N^n$ are oriented and homeomorphisms $\xi,\, \xi': S^{n-1}_{_M} \to S^{n-1}_{_N}$ reverse the natural orientations then the manifolds $M^n \#_{\xi} N^n$, $M^n \#_{\xi'} N^n$ are homeomorphic. 

  If  one of $M^n,\, N^n$ is non-orientable then for any two homeomorphisms $\xi,\, \xi': S^{n-1}_{_M} \to S^{n-1}_{_N}$ the manifolds $M^n \#_{\xi} N^n$, $M^n \#_{\xi'} N^n$ are homeomorphic. 
\end{prp}

\begin{proofn} Suppose  that  both $M^n,\, N^n$ are  orientable and   $\mu_{_M},\, \mu_{_N}$ are fixed orientations. By definition of $\xi,\, \xi'$ the composition $\xi' \xi^{-1}: S^{n-1}_{_N} \to S^{n-1}_{_N}$ is orientation-preserving. Let $C_{_N}$ be a collar of $S^{n-1}_{_N}$ in $N^n \setminus \Int B_{_N}$ and let $e: S^{n-1}_{_N} \times [0,1] \to C_{_N}$ be a homeomorphism. Hence, the homeomorphism $(e^{-1} \xi' \xi^{-1} e)|_{S^{n-1} \times \{0\}}$ preserves the orientation of $S^{n-1} \times \{0\}$. Then there is a homeomorphism $\eta: S^{n-1} \times [0,1] \to S^{n-1} \times [0,1]$ such that $\eta|_{S^{n-1} \times \{0\}} = (e^{-1} \xi' \xi^{-1} e)|_{S^{n-1} \times \{0\}}$ and $\eta|_{S^{n-1} \times \{1\}} = {\rm id}|_{S^{n-1} \times \{1\}}$. Then the formula \begin{equation}\label{homeo-id-consum}
 G([x]_{\xi}) = \begin{cases}
     [x]_{_{\xi'}},\, x \in M^{n} \setminus {\rm int}\, B^{n}_{_M};\\
     [e \eta e^{-1}(x)]_{_{\xi'}},\, x \in C_{_N};\\
     [x]_{_{\xi'}},\, x \in N^{n} \setminus {\rm int}\, C_{_N}.
 \end{cases}   
\end{equation} 
 determines the homeomorphism $G: M^n \#_{\xi} N^n \to M^n \#_{\xi'} N^n$.

Now, let us prove the proposition under the assumption that $M^{n}$ is non-orientable. If $\xi' \xi^{-1}$ preserves the orientation of $S^{n-1}_{_N}$ then the arguments are  the same as above. Suppose that  $\xi' \xi^{-1}$ reverses the orientation of $S^{n-1}_{_N}$. According to Proposition~\ref{loop}, there is the isotopy $F_t$ of $M^{n}$ such that $F_0 = {\rm id}_{_{M^{n}}}$, $F_1(S_{_M}) = S_{_M}$ and $F = F_{1}|_{S_{_M}}$ reverses the orientation of $S_{_M}$ ($S_{_N}$). By Corollary~\ref{consum-balls},  $M^{n} \#_{\xi} N^{n}$ and   $M^{n} \#_{\xi F} N^{n}$  are homeomorphic.  Hence,  it is sufficiently to prove that $M^{n} \#_{\xi F} N^{n}$  is homeomorphic to $M^{n} \#_{\xi'} N^{n}$. But $\xi' (\xi F)^{-1}$  preserves the orientation of $S^{n-1}_{_N}$ and in this case the desired homeomorphism can be constructed  similar to $G$. \end{proofn}

To complete the proof of the Theorem~\ref{group-operation}, let us recall the following classical result known as  Alexander Trick.

\begin{stat}\label{Alex}
    Let $\psi: \partial B^n \to \partial B^n$ be a homeomorphism. Then there is a homeomorphism \\ $\Psi: B^{n} \to B^n$ such that $\Psi|_{S^{n-1}} = \psi$.
\end{stat}

{\bf Proof of Theorem~\ref{group-operation}:}
The first statement of the theorem immediately follows from Corollary~\ref{consum-balls} and Proposition~\ref{consum-nonor}.  The commutativity of the connected sum operation follows from the definition. Let $B^n_{_M} \subset M^n,\, B^n_{_S} \subset S^n$ be the balls for which the connected sum is provided, $p: (M^n \setminus \Int B^n_{_M}) \sqcup (S^n \setminus \Int B^n_{_S}) \to M^n \# S^n$ be a natural projection. The map $\pi = p|_{M^n \setminus \Int B^n_{_M}}$ is a homeomorphism on the copy. Then, by Statement~\ref{Alex} $\pi$ extends to the homeomorphism $\Pi: (M^n \setminus B^n_{_M}) \sqcup (S^n \setminus B^n_{_S}) \to M^n \# S^n$.

Let us proof the associativity. Let $B^n_M\subset M^n$, $B^n_{N}, \widetilde B^n_N\subset N, B^n_{L}\subset L^n$ are locally flat balls,  $B^n_{N}\cap \widetilde B^n_{N}=\emptyset$, and $\xi_1:\partial B^n_M\to \partial B^n_N, \xi_2: \partial \widetilde B^n_N\to \partial B^n_L$ are homeomorphisms satisfying the condition of the definition of the connected sum. Then $(M^n\# N^n)\# L^n=(M^n\setminus \Int B^n_M)\cup _{\xi_1 }(N^n\setminus \Int B^n_N)\cup_{\xi_2} (L^n\setminus \Int B^n_{L})=M^n\# (N^n\# L^n)$.  Finally, it follows from  Statement~\ref{Alex}, that $M\# S^n$ is homeomorphic to $M$. \qed

The following statement shows that for some orientable manifolds  the definition of the connected sum may  be weakened.  Recall that $M, -M$ denote an oriented manifold with opposite orientations. 

\begin{prp}\label{hempel} Let $M^n, N^n$ be orientable manifolds such that at least one of $M^n,\, N^n$ admits an orientation reversing homeomorphism.  Then connected sums $M^n \# N^n$, $M^n \# (-N^n)$  are homeomorphic.   \end{prp}
\begin{proofn} Let $\mu_{_M},\, \mu_{_N}$ be fixed orientations, $B^n_{_M} \subset M^n,\, B^n_{_N}$ be locally flat  balls,   spheres  $S^{n-1}_{_M} = \partial B^n_{_M},\, S^{n-1}_{_N} = \partial B^n_{_N}$ have   orientations induced by  $\mu_{_M},\, \mu_{_N}$, respectively, and  $\xi: S^{n-1}_{_M} \to S^{n-1}_{_N}$ be an orientation-reversing homeomorphism. 

Suppose that  there exists a homeomorphism $h:M^n\to M^n$ that  reverses the orientation $\mu_{M}$.    Without loss of generality assume that  $h_{_M}(B^n_{_M}) = B^n_{_M}$. If $h_{_M}(B^n_{_M}) \neq B^n_{_M}$, then we consider a composition of $h_{_M}\widetilde h_{_M}$ instead of $h$, where  $\widetilde h_{_M}: M^n\to M^n$ is an  isotopic to identity  homeomorphism  such that $\widetilde h_{_M}(h_{_M}(B^n_{_M})) = B^n_{_M}$ (the existence of $\widetilde h_{_M}$ follows from  Proposition~\ref{balls-iso}).
Set $h_{_N}=\Id_{_N}$, and define a homeomorphism $\tilde \xi: S^{n-1}_M\to S^{n-1}_N$ by $\tilde \xi=h_{N}\xi h^{-1}_{_M}|_{S^{n-1}_{_M}}$.  Then similar to proof of Corollary~\ref{consum-balls}, one can construct a homeomorphism $g: (-M^n) \#_ {\tilde \xi} N^n \to M^n \#_\xi N^n$.  At last, we remark that  $(-M^n) \# N^n$ and  $M^n \# (-N^n)$ are homeomorphic by the  identity homeomorphism. 
Hence,  if $M^n$  admits an orientation-reversing homeomorphism then  manifolds  $M^n \# N^n$, $-M^n \# N^n$ and $M^n \# (-N^n)$ are homeomorphic. 
\end{proofn}

Any two-dimensional orientable closed manifold admit an orientation preserving homeomorphism. The Lens space $L_{3,1}$ and the complex projective plane $\mathbb{CP}^2$  are examples of 3-and 4-dimensional manifolds that do not admit such homeomorphism (see, for instance ~\cite[Lemma 3.23]{Hemp},~\cite[Exercise 1.3.1 (f)]{GoSh_en}).

Recall that an  $n$-dimensional manifold $M^n$ different from the sphere $S^n$ is called {\it prime} if it cannot be represented as a connected sum of two manifolds $M^n \# N^n$, each of which is different from the sphere $S^n$. 
In~{\cite{Mil1962} (see also~\cite[Theorem 3.22]{Hemp})} it is shown, that if  an orientable   three-dimensional manifold $M^3$ splits into connected sums $M^3=X_1 \# \cdots \# X_k=Y_1 \# \cdots \# Y_m$ of prime manifolds,  then $k=m$ and, in appropriate numeration,  $X_i$ is homeomorphic to $Y_i$ by means of an orientation preserving homeomorphism. Since $M^3 \# (-N^3)$ is homeomorphic to $(-M^3) \# N^3$, it  proves that if  $M^3 \# N^3$, $M^3 \# (-N^3)$  are  homeomorphic then at least one of $M^3, N^3$ admits an orientation reversing homeomorphism.  For $n\geq 4$ we have similar examples, in particular,  $\mathbb {CP}^2 \# \mathbb{CP}^2$ is not homeomorphic to  $(-\mathbb {CP}^2) \# \mathbb{CP}^2$ (since that manifolds have non-isomorphic intersection forms,  see~{\cite[\S 1.2]{GoSh_en}}). However,  for $n>3$ the decomposition into a connected sum  is not unique for orientable manifolds, for instance,   $\mathbb{CP}^2 \# (S^2 \times S^2), \mathbb{CP}^2 \# \mathbb{CP}^2 \# (-\mathbb{CP}^2)$ are homeomorphic while $S^2 \times S^2$, $\mathbb{CP}^2 \# (-\mathbb{CP}^2)$ are not (see~\cite{Mil1962}, \cite[Corollary 5.1.5]{GoSh_en}).

\subsection{Connected sums with $S^{n-1}$-bundles over $S^1$}\label{sec1.4}

Let  $\eta: \mathbb{S}^{n-1}\to \mathbb{S}^{n-1}$  be a homeomorphism, $n\geq 2$ and  $M^n_\eta$ be  a factor-space of $ \mathbb{S}^{n-1} \times [0,1]$ by a minimal equivalence relation such that  $(x,1)\sim (\eta(x),0)$, $x\in \mathbb{S}^{n-1}$. 
  
It is well known that  homeomorphisms  $\eta, \theta: \mathbb{S}^{n-1}\to \mathbb{S}^{n-1}$,  $n\geq 1$, are isotopic if and only if they both are either orientation-preserving or orientation-reversing (see, for instance,~\cite[\S 14.2]{GrGu-book_en} for references). Since topological type of $M^n_\eta$ depends only on isotopy class of the homeomorphism $\eta$,  there are exactly two (up to homeomorphism) manifolds of type $M^n_\eta$ and the  following statement holds (see, for instance,~\cite[Proposition~14.1]{GrGu-book_en}). 
  
\begin{stat}\label{s-gluing} If $\eta: \mathbb S^{n-1}\to \mathbb S^{n-1}$ is an orientation-preserving homeomorphism, then  $M^n_\eta$ is homeomorphic to the direct product  $\mathbb {S}^{n-1}\times \mathbb {S}^1$. If $\eta:  \mathbb S^{n-1}\to \mathbb S^{n-1}$ is an orientation-reversing homeomorphism, then $M^n_\eta$ has a structure of non-orientable $\mathbb{S}^{n-1}$-bundle over $\mathbb S^1$.    
\end{stat}

For an orientation-reversing homeomorphism $\eta$ we will denote $M^n_\eta$ by $\mathbb S^{n-1} \widetilde \times \mathbb S^1$.

The role of the manifolds $M^n_\eta$ is described in Proposition~\ref{prime-reducible}. The  purpose of this section is to prove Lemma~\ref{MeUm-nonor} and~\ref{rel-dim-n}  that are  parts of the proof of Theorem~\ref{main+}. 

The main tool of the proof is the surgery along locally flat sphere that is determined as follows.  Let  $M^n$ be a closed topological manifold  and $S^{n-1} \subset M^n$ be a locally flat sphere.  According to Statement~\ref{bicol},   $S^{n-1}$ is bi-collared,   that is there is a locally flat embedding $e:\mathbb S^{n-1} \times [-1,1]\to M^n$ such that $e(\mathbb S^{n-1}\times \{0\})=S^{n-1}$   Set $A^n=e(S^{n-1} \times [-1,1])$ and denote by  $\widehat{M}^n$  a closed manifold obtained by gluing  a copy of the  ball $\mathbb B^n$ to  each connected component of the boundary of $M^n\setminus  \Int A$. We will say that  the manifold $\widehat{M}^n$ is obtained from $M^n$ {\it by surgery along  $S^{n-1}$}. 

An {\it inverse surgery} is an operation of the obtaining a closed  manifold $M^n$ from a closed manifold $\widehat{M}^n$ by removing two disjoint locally flat balls $B^n_+, B^n_-\subset \widehat{M}^n$ and gluing the annulus $\mathbb A^n=\mathbb S^{n-1}\times [-1,1]$ to $\widehat{M}^n\setminus \Int(B^n_+\cup B^n_-)$ by means of a homeomorphism  $\theta:\partial A^n\to \partial(B^n_-\cup B^n_+)$. The result of the inverse surgery depends on an  isotopy class of $\theta$. In particular,  the following proposition holds. 

\begin{prp}\label{handle+}  Let $M^n$ be a closed manifold obtained from the sphere $\mathbb S^n$ by the inverse surgery. Then $M^n$  homeomorphic to either $\mathbb S^{n-1}\times \mathbb S^1$ or to $\mathbb S^{n-1}\widetilde \times \mathbb S^1$. 
\end{prp}
\begin{proofn}
  Let $B^n_-, B^n_+\subset \mathbb S^n$ be a locally flat balls.  By Annulus theorem, the set $\mathbb S^{n-1}\setminus \Int(B^n_-\cup B^n_+)$  is homeomorphic to $\mathbb S^{n-1}\times \mathbb S^1$. Then the proposition immediately follows from Statement~\ref{s-gluing}. \end{proofn}

\begin{figure}[!ht]
\begin{center}
\includegraphics[scale = 0.09]{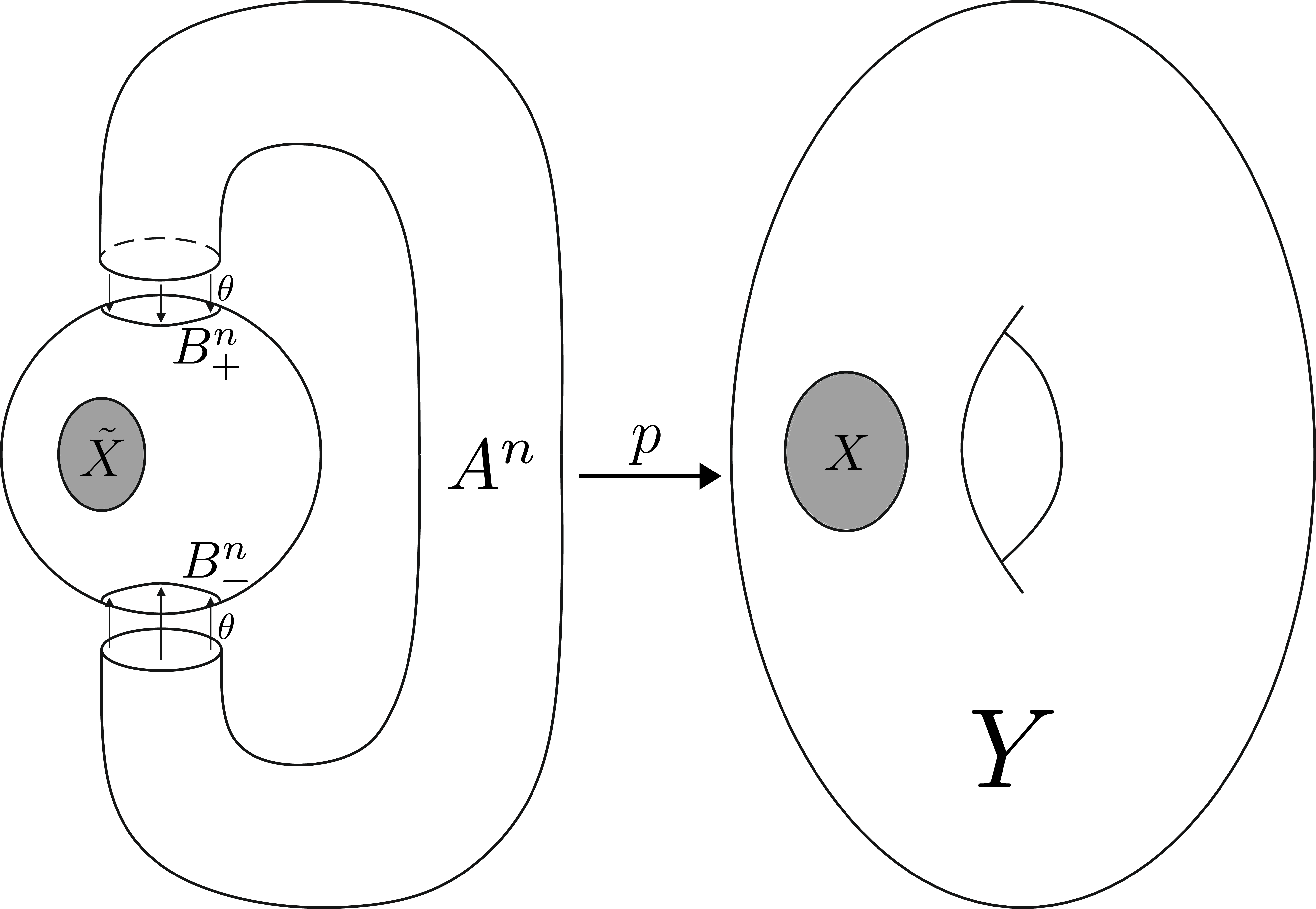}
\caption{}\label{fig-prop10}
\end{center}
\end{figure}

Due to Proposition~\ref{handle+} we  may denote  the manifold $M^n$ described there by $M_\theta$. Let $p: (S^{n-1}\setminus \Int(B^n_-\cup B^n_+))\coprod \mathbb A^n\to M_\theta$ be a natural projection,  $\widetilde X\subset S^n\setminus (B^n_-, B^n_+)$ be a locally flat ball,  $\widetilde B^n=S^n\setminus \Int X$,  and $X=p(\widetilde{X}), B^n=p(\widetilde B^n)$,  $Y=(B^n\setminus \Int(B^n_-\cup B^n_+))\cup_\theta \mathbb A^n=M_\theta\setminus \Int X$.  Since $p$ maps $\widetilde X$ into $M_\theta$ homeomorphically, then $X$  is a locally flat ball (see Fig.~\ref{fig-prop10}).   Due to Proposition~\ref{balls-iso}, we immediately get the following statement.

  \begin{cor}\label{handle-cor} Let $D^n\subset M^n_\theta$ be a locally flat ball. Then $M^n_\theta\setminus \Int D^n$ is homeomorphic to $Y$.
\end{cor}

 \begin{lem}\label{MeUm-nonor}
    Let $\widehat{M}^n$ obtained from  manifold $M^n$ by the  surgery along $S^{n-1}\subset M^n$. Then the following alternatives hold.
    \begin{enumerate}
        \item  If $\widehat{M}^n$ is disconnected, then $M^n=M_1 \# M_2$ where $M_1, M_2$ are connected components of $\widehat{M}^n$;
        \item if $\widehat{M}^n$ is connected and $M^n$  is orientable, then $M^n$ is homeomorphic to  $\widehat{M}^n \# (S^{n-1} \times S^1)$;
        \item if $\widehat{M}^n$ is connected and $M^n$  is non-orientable then $M^n$ is homeomorphic either to $\widehat{M}^n \# (S^{n-1}  \times S^1)$ or   to     $\widehat{M}^n \# (S^{n-1} \widetilde \times S^1)$.
    \end{enumerate}
\end{lem}
\begin{proofn}
If $\widehat{M}^n$ is disconnected, then the statement immediately follows from the definition of the connected sum. If $\widehat{M}^n$ is connected and $M^n$ is oriented,  then the statement  is proved in~\cite[Lemma 7]{MeUm79}. So, we have to prove  item 3.  Suppose that $M^n$ is non-orientable  and $\widehat{M}^n$ is connected. 

By definition, there are disjoint  locally flat  balls  $B^n_-,\, B^n_+\subset \widehat {M}^n$ such that $M^n$ is the result of gluing the annulus  $\mathbb A^n=\mathbb S^{n-1}\times [-1,1]$ to $\widehat{M}^n \setminus \Int (B^n_-\cup B^n_+)$ by means of a homeomorphism  $\theta: \partial \mathbb A^n \to \partial\, (B^n_-\cup B^n_+)$. Let $C\subset \widetilde M^n\setminus \Int B^n_-$ be a collar of $\partial B^n_-$. It follows from Proposition~\ref{balls-iso} that there exists a homeomorphism $h:\widehat{M}^n\to \widehat{M}^n$ such that $h|_{_{B^n_-}}=\Id_{_{B^n_-}}$ and $h(B^n_+)\subset \Int C$. So, without loss of generality we may assume that $B^n_+\subset \Int C$ (otherwise consider $h(B^n_+)$ and $h\theta$ instead of $B^n_+$ and $\theta$). Hence, $B^n_-\cup B^n_+$ belongs to the interior of a ball $B^n=B^n_-\cup C$. A manifold $(B^n\setminus \Int(B^n_-\cup B^n_+))\cup_\theta \mathbb A_n$ is homeomorphic to the manifold $Y=M^n_\theta\setminus D^n$ described in Corollary~\ref{handle-cor}. Let $\xi: (B^n\setminus \Int(B^n_-\cup B^n_+))\cup_\theta \mathbb A_n\to Y$ be a homeomorphism, and  $i:\partial B^n\to \widehat M$ is an inclusion map.   Then $M^n=(\widehat M^n\setminus \Int(B^n_-\cup B^n_+))\cup_\theta \mathbb A^n=(\widehat M^n\setminus \Int\,B^n)\cup_i (B^n\setminus \Int(B^n_-\cup B^n_+))\cup_\theta \mathbb A^n=(\widehat M^n\setminus \Int\,B^n)\cup _\xi (M^n_\theta\setminus D^n)=M^n\# M^n_\theta$. \end{proofn}

The following proposition is the generalisation of well-known facts in low dimension (see, for instance,~\cite[Lemma 3.17, Theorem 3.21]{Hemp}). 

\begin{prp}\label{rel-dim-n}
    Let $M^n$ be a  non-orientable  closed manifold of dimension  $n>3$. Then  \begin{equation}
    M^n \# (S^{n-1} \times S^1) \cong M^n \# (S^{n-1} \widetilde \times S^1).
    \end{equation}
\end{prp}
\begin{proofn}
Recall that $M^n_\eta$ denotes either $S^{n-1} \times S^1$ or $S^{n-1} \widetilde \times S^1$. Set $Q^n=M^n\# M^n_\eta$.  
Due to Theorem~\ref{group-operation}, $M^n\# M^n_\eta=M^n\# S^n\# M^n_\eta$. Hence, manifold $Q^n$ can be considered as a result of the inverse surgery on $M^n$. That means, that there are   locally flat balls $B^n_-,\, B^n_+\subset M^n$   and  a homeomorphism  $\theta: \partial \mathbb A^n \to \partial\, (B^n_-\cup B^n_+)$ such that $Q^n$ is homeomorphic to $(M^n\setminus \Int (B^n_-\cup B^n_+))\cup _\theta \mathbb A^n$.

 There are two possibilities: $M^n_\eta$ is either  orientable or not,  that depends  only on the   isotopy class of gluing map $\theta$.   It follows from  Proposition~\ref{loop} that  there exists a  homeomorphism $F: M^n\setminus \Int(B^n_-\cup B^n_+)$ such that  $F|_{\partial{B^n_-}}$ is identity and $F|_{\partial{B^n_+}}$ is the orientation reversing.  Set $\widetilde \theta=F\theta$. Then $\theta$ and $\widetilde \theta$ are non-isotopic and exhaust all possible isotopic classes of sphere homeomorphisms and all possible topological types of $Q^n$. 
 
 Determine a map  $G: (M^n \setminus \Int (B^n_0 \cup B^n_1)) \cup_\theta \mathbb A^n\to (M^n \setminus \Int (B^n_0 \cup B^n_1)) \cup_{\widetilde\theta} \mathbb A^n$ by  \begin{equation}
 G([x]_\theta) = \begin{cases}
     [x]_{\widetilde\theta},\, x \in \mathbb A^n;\\
     [F(x)]_{\widetilde\theta},\, x\in M^n\setminus \Int(B^n_0\cup B^n_1).
 \end{cases}
 \end{equation}

 If $x\in \partial \mathbb A_n$, then  $[x]_{ \theta}=\{x, \theta(x)\}$  and $[x]_{\widetilde \theta}=\{x, \widetilde \theta(x)\}=\{x, F(\theta (x)\}=G([x]_\theta)$. Hence, $G$ is a homeomorphism. 
 \end{proofn}

A manifold $M^n$ is said to be {\it irreducible} if any locally flat sphere $S^{n-1}\subset M^n$ bounds a ball $B^n\subset M^n$. Any prime and non-irreducible three-dimensional manifold is homeomorphic to $M^3_\eta$ (see the proof in~\cite[Lemma 3.8]{Hemp}). We prove the following generalisation.

 \begin{prp}\label{prime-reducible}
     Any prime and non-irreducible manifold of dimension $n> 3$ is homeomorphic to $M^n_\eta$.  
 \end{prp}
 \begin{proofn}
     Suppose that  $M^n$ is prime and non-irreducible. Then  any locally flat $(n-1)$-dimensional sphere in $M^n$ either bounds a ball in $M^n$ or does not separate $M^n$. Let $S^{n-1} \subset M^n$ be a locally flat sphere that does not separate $M^n$.  Then a manifold $\widehat M^n$ obtained from $M^n$ by the surgery along $S^{n-1}$ is connected. Due to Lemma~\ref{MeUm-nonor}, 
     $M^n=\widehat M^n\# M^n_\eta$. But $M^n$ is prime, hence, by Theorem~\ref{group-operation}, $\widehat M^n$ is homeomorphic to $S^{n-1}$ and $M^n$ is homeomorphic to $M^n_\eta$. 
     
 \end{proofn}

 \section{Regular homeomorphisms
}\label{sec1.5}

Let $f: M^n \to M^n$ be a homeomorphism on a closed topological manifold $M^n$ with a metric $\rho$. A point $x, y\in M^n$  is said to be connected by  {\it  an $\varepsilon$-chain of $f$}  if there   is a sequence of points  $x = x_0, \dots , x_k = y$  and a  sequence $m_1, \dots, m_k$ of integers  such that $\rho(f^{m_i}(x_{i-1}), x_i) < \varepsilon$, $m_i \geq 1$ for $1 \leq i \leq k$. A number  $N=m_1 + \dots + m_k$ is called a length of the chain.

A point $x \in M^n$ is {\it chain recurrent} for the homeomorphism $f$ if for any $\varepsilon > 0$ there are $N$ and $\varepsilon$-chain of length $N$ connecting $x$ to itself. The set $\mathcal R_f$ of all chain recurrent points of $f$ is the {\it chain recurrent set}, and  its connected components are {\it chain components}. It immediately follows form the definition, that the chain recurrent set is $f$-invariant, hence, the following statement holds.   

\begin{prp}\label{ch-recc-hyperbolicity} If the chain recurrent set of a homeomorphism  $f:M^n\to M^n$ is finite, than it consists of periodic points. 
    \end{prp}

Suppose that  $M^n$ is smooth,  $f:M^n\to M^n$ is a diffeomorphism, and $p$ is its periodic point of  period $m_p$. The point $p$  is called {\it hyperbolic}  the   differential $D_pf^{m_p}: \mathbb R^n \to \mathbb R^n$ of $f^{m_p}$ have no eigenvalues with absolute value equal to one. A number $i_p\in \{0,\dots, n\}$ of eigenvalues of $D_{p}f^{m_p}$  with absolute value greater that one is called  {\it a Morse index of $p$}. It follows from Grobman-Hartman Theorem {(see~\cite{Gr59, Gr62, Hart63}))} and a classification of linear automorphisms of $\mathbb R^n$ {(see~\cite[Proposition 2.9]{PaMe}))}, that 
\item[($*$)] there exists a  
neighbourhood $U_p \subset M^n$ of $p$  and a homeomorphism $h_p: U_p\to \mathbb R^n$ such that  $f^{m_p}|_{U_p}=h^{-1}_pa_{i_p, {\delta_u, \delta_s}}h_p|_{U_p}$, where   $a_{i_p, {\delta_u, \delta_s}}:\mathbb R^n\to \mathbb R^n$  is a map defined by \begin{equation}
a_{i_p, {\delta_u, \delta_s}}(x_1,..., x_{i_p}, x_{{i_p}+1},..., x_n) = \left(2 \delta_u x_1,2x_2,..., 2x_{i_p}, \dfrac{\delta_s}{2}x_{{i_p}+1},\dfrac{1}{2}x_{i_p+2}, ..., \dfrac{1}{2}x_{n}\right)
\end{equation}
and $\delta_u,\delta_s \in \{+1,-1\}$.

We call a periodic point $p$ of the homeomorphism $f$  {\it topologically hyperbolic} if the condition~($*$) holds. 
  The number $i_p$ is called  {\it a Morse index of $p$}. If $i_p=0$ then $p$ is called {\it a sink}, if $i_p=n$ then $p$ is called {\it a source}, otherwise $p$ is called {\it a saddle} periodic point. 

The map  $a_{i_p, \delta_u, \delta_s}$ induces a splitting  $\mathbb R^n=E^u_{i_p} \oplus E^s_{i_p}$ on invariant  linear subspaces $E^u_{i_p},\, E^s_{i_p}$ of dimensions $i_p,\, (n - i_p)$, respectively, that are unstable and stable manifolds of a fixed point $O$.  A set $W^u_{p, loc} = h^{-1}(E^u_{i_p})$ is called {\it a local unstable  manifold}, and a set  $W^u_p = \bigcup\limits_{i \in \mathbb Z} f^i(W^u_{p, loc})$  is called {\it an unstable manifold} of the topologically hyperbolic periodic point $p$. A local and global stable manifolds    $W^s_{p, loc}, W^s_{p}$  of the topologically hyperbolic periodic point $p$ is defined as the local and global unstable manifolds of $p$ with respect to $f^{-1}$. A connected component $\ell^u_p\, (\ell^s_p)$ of the set $W^u_p\setminus p\, (W^u_p\setminus p)$ is called {\it an unstable $($stable$)$ separatrix} of $p$.

Due to~\cite[Statement~1]{PoZi_RCD},
\begin{equation}
W^u_p=\{q\in M^n|\,  \lim \limits_{n\to +\infty} \rho(f^{-n m_p}(q),p)=0\}, \\  W^s_p=\{q\in M^n|\, \lim\limits_{n\to +\infty}\rho(f^{n m_p}(q),p)=0\}.
\end{equation}

A  homeomorphism $f: M^n\to M^n$  is called {\it regular} if its chain recurrent set $\mathcal R_f$ is finite and topologically hyperbolic. 

Morse-Smale diffeomorphism and gradient-like flows are important and motivating example of the  regular dynamical systems. In general, trajectories of regular  dynamical systems have more complex asymptotic  behavior than ones for Morse-Smale systems  since we omit a requirements of the  transversality of the intersection of invariant manifolds. The following statements describe properties of regular homeomorphisms that they share with  Morse-Smale diffeomorphisms. First two are proved in~\cite[Statement 2]{PoZi_RCD},~\cite[Theorem 1]{PoZi_RCD}  for regular homeomorphism and in~\cite[Statement 1.2.5]{GrPo-book}, \cite[Theorem 2.3]{Sm}, for Morse-Smale diffeomorphism. 

Following~Smale, we determine {\it a Smale relation} $\prec$ on the set $\mathcal R_f$ by the rule: $q \prec p$ if and only if $W^u_p \cap W^s_q \not= \emptyset$. 

\begin{stat}\label{reghom-order} 
Let $f: M^n \to M^n$ be a regular homeomorphism. Then 
\begin{enumerate}
    \item for any  points $p,\, q,\, r\in \mathcal {R}_f$ conditions $p \prec q,\, q \prec r$ imply $p \prec r$;
    \item there is no set of pairwise distinct points  $p_1, \dots, p_k\in \mathcal {R}_f$ such that  $p_i \prec p_{i+1}$ for  any $i \in \{1, \dots, k-1\}$  and $p_k \prec p_1$.
\end{enumerate}
\end{stat}

\begin{stat}\label{reghom-smale}
    Suppose $f: M^n \to M^n$ to be a regular homeomorphism. Then:
    \begin{enumerate}
        \item $M^n= \bigcup\limits_{p \in \mathcal R_f}W^u_p = \bigcup\limits_{p \in \mathcal R_f}W^s_p$; 
        \item for any periodic point $p \in \mathcal R_f$ of Morse index $i_p$ the set $W^u_p$ is topological submanifolds of $M^n$ homeomorphic to $\mathbb R^{i_p}$;
        \item $({\rm  cl}\, W^u_p) \setminus W^u_{p} \subset  \bigcup\limits_{q\prec p}  W^u_q$ for any $p \in \mathcal{R}_f$.  
    \end{enumerate}
\end{stat}

{\begin{cor}\label{tube} Let $p\in M^n$ be a topologically hyperbolic point of period $m_p$,  $B^{i_p}_p\subset W^u_p$ be a  compact ball such that $p\in \Int B^{i_p}_p$. Then there is a  compact neighborhood $V_p\subset M^n$ of $B^{i_p}_p$  and a homeomorphism  $g:\mathbb B^{i_p}\times \mathbb B^{n-i_p}\to V_p$ such that $g(\mathbb B^{i_p}\times\{O\})=B^{i_p}_p$, $g(\{O\}\times \mathbb B^{n-i_p})=V_p\cap W^s_p$, and  projections $\pi_u: V_p\to W^{u}_p, \pi_s: V_p\to W^s_p$ along the fibers have  the following properties:
\begin{enumerate}
\item $\pi_u^{-1}(p)\subset W^s_p$,  $\pi_s^{-1}(p)=B^{i_p}_p$;
\item $f^{-k m_p}(\pi_{u}^{-1}(x))\subset \pi_{u}^{-1}(f^{-k m_p}(x))$, $f^{k m_p}(\pi_{s}^{-1}(y))\supset \pi_{s}^{-1}(f^{k m_p}(y))$ for any $x\in B^{i_p}_p,  y\in V_p\cap W^s_p$ and  $k\in \mathbb N$. 
\end{enumerate} 
 \end{cor}}
\begin{proofn}
Let $U_p, h_p:U_p\to \mathbb R^n$ be  a neighborhood of $p$ and a homeomorphism satisfying condition~$(*)$.  Remark that for any balls $D^u\subset E^u_{i_p}, D^s\subset E^s_{i_p}$ of dimension $i_p, (n-i_p)$, respectively,  containing the origin $O$, a set $V=D^u\times  D^s$  has the  properties described in the statement of the corollary, with the  formal replacing $f^{m_p}$  with $a_{i_p, \delta_u, \delta_s}$. 
Since $B^{i_p}\subset W^u_{p}$, there is $\nu$ such that $f^{-\nu m_p}(B^{i_p})\subset U_p$.     Set $D^u=h_p(f^{-\nu m_p}B^{i_p})$. Due to Statement~\ref{reghom-smale}, Statement~\ref{reghom-order},  $W^u_{p}$ is the submanifold of $M^n$  and $W^u_p\cap W^s_p=p$. Hence   we may choose  $D^s$ in such a way that $h_p^{-1}(V)\cap (W^u_{p}\cup W^s_p)=h_p^{-1}(D^u\cup D^s)$. Then  the neighbourhood $h_p^{-1}(V)$  has the fibre structure with the required  properties, and so does $V_p=f^{-\nu m_p}(h_p^{-1}(V))$.  
\end{proofn}

Set $W^s_{_{\mathcal O_l}}=\bigcup \limits_{p\in \mathcal O_l}W^s_p$, $W^u_{_{\mathcal O_l}}=\bigcup \limits_{p\in \mathcal O_l}W^u_p$. 
We say that $\mathcal O_p\prec \mathcal O_q$  if  $x \prec y$ for some $x\in \mathcal{O}_p, y\in \mathcal O_q$. Due to Statement~\ref{reghom-order}, Smale relation can be extended to a total order relation $\preccurlyeq$ on the set of periodic  orbits $\{\mathcal O_p, p\in \mathcal R_f\}$ of the regular homeomorphism  $f$ as follows: $\mathcal O_p\preccurlyeq \mathcal O_q$  if and only if either $\mathcal O_p\prec \mathcal O_q$ or $W^s_{\mathcal O_p}\cap W^u_{\mathcal O_q}=\emptyset$. 

Suppose that  all periodic orbits consisting in $\mathcal R_f$ are numbered  with respect to the total order: \begin{equation}\label{total-order}
 \mathcal O_{1} \preccurlyeq   \dots \preccurlyeq  \mathcal O_{k}.     
 \end{equation}

It follows from Statement~\ref{reghom-smale} that $\mathcal R_f$ contains at least one sink and one source periodic orbits. Then we may assume that first $\mu$ orbits in the row~(\ref{total-order})  are orbits of all sink periodic points from $\mathcal R_f$.   Below we suppose  that   a set of saddle orbits of $f$ is non-empty,  otherwise $\mathcal R_f$ consists of exactly one sink and one source fixed points and $M^n$ is homeomorphic to the sphere $S^n$ as it shown in~Proposition~\ref{polar}.
 
Set $A_j=\bigcup \limits_{l=1}^j W^u_{_{\mathcal{O}_l}}, R_j=\bigcup \limits_{l=j+1}^{k} W^s_{_{\mathcal{O}_l}}$, $j\in \{1,\dots, k-1\}$.  We  call a maximal dimension of the unstable (stable) manifolds of the orbits $\mathcal O_l\in A_j\ (R_j)$   the {\it dimension of the set $A_j \ (R_j)$}.

Recall that a set $A\subset M^n$ is called {\it an attractor} of $f$ if there exists a compact neighbourhood ({\it a trapping neighbourhood}) $Q_a\subset M^n$ of $A$ such that $f(Q_a)\subset \Int Q_a$ and $A= \bigcap\limits_{n\in \mathbb Z}f^n(Q_a)$.
 A set $R\subset M^n$ is {\it a repeller} of $f$ if it is an attractor for $f^{-1}$. 
 
The following proposition is similar to~\cite[Theorem 1.1]{GrZhMePo2010}, but is proved bypassing the smooth technique.

\begin{theo}\label{trap} $A_j$ is an  attractor  with a trapping neighbourhood $Q_j\subset \bigcup \limits_{l=1}^{j} W^s_{_{\mathcal O_l}}$ such that 
$\partial Q_j$ is a locally flat submanifold of $M^n$. If $j\geq \mu$,  for any 
 $p\in \mathcal O_{j+1}$   there exists a compact ball $B^{i_p}_p\subset W^u_p$   such that $p\in \Int B^{i_p}_p$,  $W^u_p\setminus B^{i_p}_p\subset \Int Q_j$ and $\partial Q_j$ is  embedded topologically transversely to $B^{i_p}_p$. 
\end{theo}
\begin{proofn}
We construct the trapping neighbourhood for $A_j$ using induction by $j\subset \{1,\dots, k-1\}$.  

Let $j=1$,  $\omega\subset \mathcal O_1$ be a  sink periodic  point  of period  $m_\omega$, and   $U_\omega, h_\omega:U_\omega\to \mathbb R^n$ be a neighbourhood of $\omega$ and a homeomorphism satisfying the condition $(*)$. Set $B^n_\omega=h_\omega^{-1}(\mathbb B^n)$. Since  $h_\omega$ conjugates $f^{m_\omega}|_{_{U_\omega}}$ with $a_{0, \delta_s, \delta_u}$ and $a_{0,\delta_s, \delta_u}(\mathbb B^n)\subset \Int \mathbb B^n$, then $f^{m_p}(B^n_\omega)\subset \Int B^n_\omega$.  

Let $p\in \mathcal O_{\mu+1}$ be a saddle periodic point. Then  $W^u_p\setminus p\subset \bigcup\limits_{l=1}^\mu W^s_{\mathcal O_l}$. Let $B^{i_p}_p\subset W^u_p$ be a compact ball such that $p\in \Int B^{i_p}_p$. Due  to Corollary~\ref{tube}, $B^{i_p}_p$ has a trivial microbundle in $M^n$. Denote by  $\ell^u_p$ a connected component of the set $W^u_p\setminus p$. Since $W^s_{\omega}\cap W^s_{\omega'}=\emptyset$ for a any sink periodic point $\omega'\neq \omega$, there exists a single sink $\omega$ such that $\ell^u_p\subset W^s_\omega$.  There are two cases: 1)  $i_p>1$, 2) $i_p=1$. In case 1) the set $W^{u}_p\setminus p$ is connected and hence $\partial B^{i_p}_p\subset W^s_\omega$. Then  there is a number $N$ such that $f^{n m_p}(\partial B^{i_p}_p)\subset \Int B^n_\omega$ for all $n>N$, and, consequently,   $W^u_p\cap \partial B^n_\omega=B^{i_p}_p\cap \partial B^n_\omega$.  According to  Statement~\ref{quinn}, there is an isotopy $h_t:M^n\to M^n$   such that the ball $Q_{\omega}=h_1(B^n_\omega)$ is embedded topologically transversely to $B^{i_p}_p$ and $f^{m_\omega}(Q_\omega)\subset \Int Q_\omega$.

In case 2),  $\partial B^{i_p}_p$ consists of two points ${z_+, z_-}$ that lays on different connected components of $W^u_p\setminus p$,  and, possible, in stable manifolds of different sink periodic points $\omega_+, \omega_-$. If $\omega_+, \omega_-$ belongs to different orbits,  we choose the ball $Q_{\omega_{\pm}}\subset W^s_{\omega_{\pm}}$  similar to case  1).  If there exists a number $m$ such that $f^m(\omega_-)=\omega_+$, we choose a ball $Q_{\omega_-}$ and set $Q_{\omega_+}=f^m(Q_{\omega_-})$.

For any orbit $\mathcal O_l$, $l\in \{1,\dots, \mu\}$, choose a point $\omega\subset \mathcal O_l$ and set 
 $Q_{l}=\bigcup \limits_{m=0}^{m_{\omega}-1}f^{m}(Q_{\omega})$. 
By construction,   $Q_l$ is the desired trapping neighbourhood for the set $A_l$, $l\in \{1,\dots, \mu\}$.

Suppose that we built a trapping neighborhood $Q_{j}$ for $A_j$, $j\subset \{\mu,\dots, k-2\}$. Let us construct the trapping neighbourhood for $A_{j+1}$.  There are two possibilities: 1) $\mathcal O_{j+1}$ is the orbit of a source periodic  point; 2) $\mathcal O_{j+1}$ is the orbit of a saddle periodic point. 

In case 1) set $Q_{j+1}=Q_j\cup W^u_{_{\mathcal O_{j+1}}}$.  Consider case 2). Let $\sigma\in \mathcal O_{j+1}$ be  a saddle periodic point with period $m_\sigma$ and Morse index $i_\sigma$;  $B^{i_\sigma}_\sigma, V_\sigma$  be  a compact ball and its compact neighbourhood satisfying the conclusion of Corollary~\ref{tube}.  By Corollary~\ref{tube}, there exists a homeomorphism $g:  \mathbb B^{i_\sigma}\times \mathbb B^{n-i_\sigma}\to V_\sigma$ such that  $g(B^{i_\sigma}\times \{O\})=B^{i_\sigma}_\sigma$. Set $T_1=g(\partial \mathbb B^{i_\sigma}\times \mathbb B^{n-i_\sigma}), T_2=g(\mathbb B^{i_\sigma}\times \partial \mathbb B^{n-i_\sigma})$. Then $\partial V_\sigma=T_1\cup T_2$. Remark that by Corollary~\ref{tube} $f^{m_p}(T_2)\cap V_\sigma$ has a trivial microbundle in $M^n$. Then so does $T_2$.  

\begin{figure}[!ht]
\begin{center}
\includegraphics[scale = 0.15]{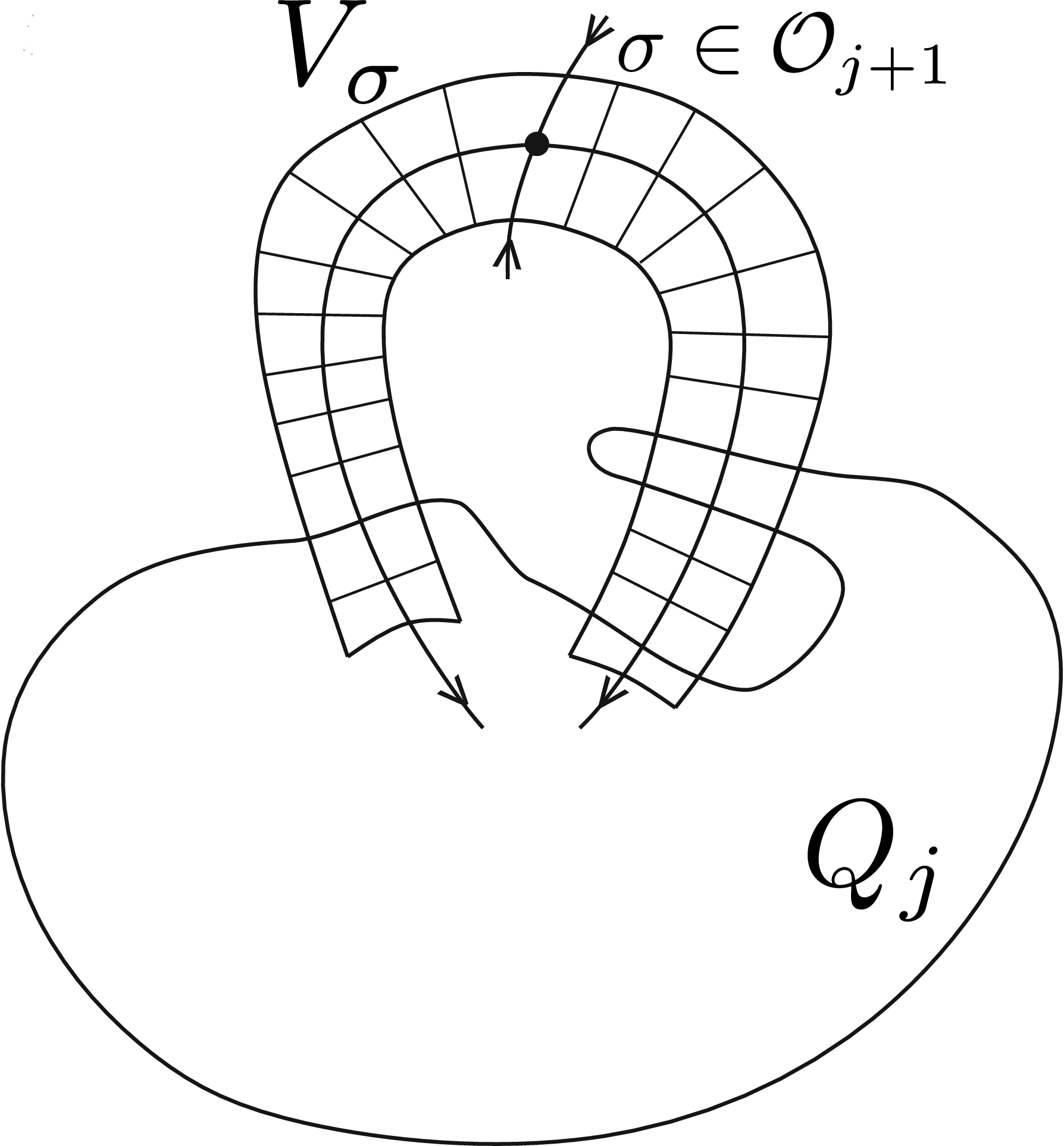}
\caption{A construction of the trapping neighborhood for $A_{j+1}$}\label{tr}
\end{center}
\end{figure}

Since $\sigma\in \mathcal O_{j+1}$, $W^u_\sigma\subset \bigcup \limits_{l=1}^j  W^s_{\mathcal O_l}$ and, moving to an iteration $f^N(B^{i_\sigma}_\sigma)$, if necessary, we may suppose that $T_1\subset \Int Q_j$. Applying Statement~\ref{quinn} once more, we  change $\partial Q_j$ near a neighbourhood of $\partial Q_j\cap T_2$ by an ambient isotopy $\hat h_t: M^n\to M^n$ so that $\hat h_1(\partial Q_j)$   is embedded topologically transversal to $T_2$. We keep a notation $Q_j$ for $\hat {h}_1(Q_j)$.  
Then  $X=\partial Q_{j}\cap T_2$ is a submanifold of $\partial Q_j$ and $T_2$. Moreover, $X$  is a boundary of $\partial Q_i\setminus \Int V_\sigma$ and of $\partial V_\sigma\setminus \Int Q_j$. Then the set  $Q_{\sigma}=Q_j\cup V_\sigma$ is bounded by a locally flat manifold $(\partial Q_i\setminus \Int V_\sigma)\cup (\partial V_\sigma\setminus \Int Q_j)$ (see Fig.~\ref{tr}).

It follows from Corollary~\ref{tube} that $f^{m_\sigma}(V_\sigma\setminus Q_j)\subset \Int(V_\sigma\cup Q_j)$. 

Set $V_{_{\mathcal O_{j+1}}}=\bigcup \limits_{m=0}^{m_\sigma-1}f^m(V_\sigma)$ and $Q_{j+1}=Q_j\cup V_{_{\mathcal O_{j+1}}}$. By construction, $Q_{j+1}$ is the desired trapping neighborhood for $A_{j+1}$.

\end{proofn}
\begin{cor}\label{dim} $R_j$ is a repeller of the regular homeomorphism $f: M^n\to M^n$.  If $\Dim R_j\leq (n-2)$ then $A_j$ is connected.
\end{cor}
To prove  that $R_j$ is the repeller it is enough to  apply the Theorem~\ref{trap}   to  $f^{-1}$.  The proof of the  connectivity of $A_j$ for  appropriate dimension of $R_j$  is literally the same as in~\cite[Theorem 1.1]{GrZhMePo2010}. The idea of the proof is the following. Since $\Dim R_j<n-1$, it does not divide $M^n$. Then $\bigcup \limits_{l=1}^{j} W^s_{_{\mathcal O_l}}=M^n\setminus R_j$ is connected and so does $Q_j$. Hence, $A_j$ is connected as an intersection of nested connected sets $\{f^m(Q_j)\}$.  

\section{Proof of Theorem~\ref{main+}}\label{tm1}
\subsection{Polar regular homeomorphisms}

Recall that $k_i$ denotes a number of periodic points of the regular homeomorphism $f: M^n\to M^n$ whose Morse index equals to $i\in \{0, \dots, n\}$.  
In this section we provide a proof of~Theorem~\ref{main+} for  case $k_1=k_{n-1}=0$.
 For smooth systems it follows from~\cite[Corollary 1.1]{GrZhMePo2010} and~\cite[Theorem 1.3]{GrGuPo15}\  
 (see also~\cite[Proposition 4.1]{GrGu-book_en}).    

\begin{prp}\label{polar} Let $M^n$  be a connected closed  topological   manifold of dimension $n\geq 2$ and $f: M^n\to M^n$ be  a regular homeomorphism  whose  non-wandering set does not contain periodic  points with Morse index equal to $1$ and  $(n-1)$. Then 

\begin{enumerate}
    \item $f$ is polar;
    \item if $f$ has no heteroclinic intersections, then $M^n$  is homeomorphic to a sphere $S^n$ if and only if the set of saddle points of $f$ is empty;
    \item $M^n$ is simply connected. 
\end{enumerate}
\end{prp}
\begin{proofn}
Let us prove Item 1. Let $A_\mu$  be a union of all sink orbits of $f$. Set  $R_\mu=\bigcup \limits_{p\in \mathcal R_f\setminus A_\mu}W^s_p$. By the conditions,  $\Dim R_\mu<n-1$. Then due to Corollary~\ref{dim}, $A_\mu$  is connected, and, consequently, consists of a single point $\omega$. Similar arguments for $f^{-1}$ prove that the set of all source periodic orbits of $f$ is also consists of a single point $\alpha$. Hence, $f$ is polar. 

The proof of Item 2 is completely similar to the proof of \cite[Theorem 1.3]{GrGuPo15}. We repeat here the main idea of it. If  $f$ has no any  saddle points, then its non-wandering set contains exactly one source $\alpha$ and one sink $\omega$ and  $M^n=W^u_\alpha\cup \{\omega\}$. Since $W^u_\alpha$ is an  open ball of dimension $n$, $M^n$ is homeomorphic to the sphere $S^n$.  Suppose that  $M^n$ is the sphere, the set of saddle periodic points of $f$ is non-empty and invariant manifolds of any saddle periodic point $p$ does not intersect invariant manifolds of other saddles. Then, due to Statement~\ref{reghom-smale}, $ \Cl W^u_p=W^u_p\cup \omega, \Cl W^s_p=W^s_p\cup \alpha, $. Hence, $\Cl W^u_p, \Cl W^s_p$ are spheres of dimension $i_p, n-i_p$, respectively, that intersect each other  at a single point $p$. Then the intersection number of $\Cl W^u_p, \Cl W^s_p$ is different from zero. On the  other hand,  there is a sphere $S^{i_p}\subset M^n$ such that $\Cl W^s_p\cap S^{i_p}=\emptyset$, and, consequently, the intersection number  of $\Cl W^s_p,  S^{i_p}$ equals zero.  Since $H_{i_p}(S^n)=0$ for $i_p\in \{2, \dots, n-2\}$,  $S^{i_p}$ is homological to $\Cl W^u_p$. Since  the intersection number is a homological invariant,  we obtain the contradiction that proves that the set of saddle periodic points  of $f$ in this case is empty.

Let is prove Item 3. 
Suppose that  the set of saddle periodic points of $f$ is not empty and consists of points of Morse indices $2,\dots, (n-2)$. Then $n\geq 4$. Let  $\gamma \subset M^n$ be  a loop representing a class $[\gamma]\in \pi_1(M^n)$.  It follows from~\cite{GenPos}, that $\gamma$ may be considered as a locally flat embedded circle. Let us show that $\gamma$ can be moved by an ambient isotopy of $M^n$ to a loop $\tilde\gamma\subset W^u_\alpha$. Since $W^u_\alpha$ is homeomorphic to $\mathbb R^n$ and, consequently, is simply connected,  $\tilde \gamma$ is homotopic to zero, and so would be  $\gamma$, that meant  $\pi_1(M^n)$ is trivial. 

We will construct the  desired ambient isotopy moving $\gamma$ sequentially outward unstable manifolds of all orbits $ \omega=\mathcal O_1\preccurlyeq\dots \preccurlyeq \mathcal O_{k-1}$ preceding $\alpha$.

Due  to Statement~\ref{quinn}, there is an isotopy $h_t: M^n\to M^n$ such that $h_1(\gamma)\cap \omega=\emptyset$. Then  there is a trapping neighbourhood $Q_1$  of $\omega$  such that $h_1(\gamma)\cap Q_1=\emptyset$. Set $\gamma_1=h_1(\gamma)$. 
It follows from Theorem~\ref{trap} that  $X=W^u_{\mathcal O_2}\cap (M^n\setminus \Int Q_1)$ is compact and belongs to  a union of compact balls laying  in $W^u_{\mathcal O_2}$. Then, by Corollary~\ref{tube}, $X$ has a trivial normal microbundle in $M^n$. Using Statement~\ref{quinn} once more, we construct an isotopy $g_t: M^n\to M^n$ relative $Q_1$ and such that $g_1(\gamma_1)$  is embedded topologically transversely to  $X$. Since $\Dim X+\Dim \gamma_1\leq n$, it means that $X\cap g_1(\gamma_1)=\emptyset$. Hence, $g_1(\gamma_1)\cap W^u_{\mathcal O_2}=\emptyset$. Due to Theorem~\ref{trap}, $\omega\cup W^u_{\mathcal O_2}$  is the attractor of $f$ with a trapping neighborhood $Q_2$. For a sufficiently large $N$, we have $f^{N}(Q_2)\cap g_1(\gamma_1)=\emptyset$. Repeating the arguments above, after a finite number of applying Statement~\ref{quinn} and Theorem~\ref{trap}, we get the desired isotopy. 
\end{proofn}

\subsection{A surgery along $(n-1)$-dimensional separatrices}\label{sec1.6}

Previous section proves Theorem~\ref{main+} for the  case $k_1=k_{n-1}=0$.  If $k_1^2+k_{n-1}^2\neq 0$, we will use the surgery along locally flat sphere, introduced in Section~\ref{sec1.4}. The following statements show that a closure of codimension one    separatrix is a suitable sphere, more over, we may determine a regular homeomorphism from class $G(\widehat{M}^n)$  on the resulting manifold $\widehat{M}^n$. 

Let $\sigma$ be a saddle periodic point of a regular homeomorphism $f: M^n\to M^n$ such that $\Dim W^u_{\sigma}=n-1$  and $W^u_{\sigma}$ does not contain heteroclinic points. It follows from Statement~\ref{reghom-smale} that  there is a unique sink periodic  point $\omega$ such that  $({\rm cl}\, W^u_{\sigma})\setminus \sigma\subset W^s_\omega$ and  ${\rm cl}\, W^u_{\sigma}=W^u_{\sigma}\cup \omega$. Set $S_{\sigma} = W^u_{\sigma_{n-1}} \cup \omega$.

\begin{lem}\label{locfl-sep}
    $S_{\sigma}$  is  a locally flat $(n-1)$-dimensional sphere.  
\end{lem}
\begin{proofn} Due to Statement~\ref{reghom-smale}, $W^u_{\sigma}$ is a submanifold of $M^n$ homeomorphic to $\mathbb R^n$.  Then ${\rm cl}\, W^u_{\sigma}$ is a sphere of dimension $(n-1)$ which is locally flat at all points except, possibly,  $\omega$. It follows from~\cite[Theorem 1]{Kir68} (see also \cite[Corollary 3A.6]{Dav})   that for $n \geq 4$ the sphere $S^{n-1}\subset M^n$ cannot have one point of wildness (in fact, the set of points of wildness no less than  uncountable). Then $S_{\sigma}$ is locally flat at $\omega$. 
\end{proofn}

\begin{cor}\label{cilnbh}  There is a neighbourhood $N_{\sigma}\subset W^s_\omega\cup W^s_{\sigma}$ of $S^{\sigma}$ homeomorphic to $S^{n-1}\times [-1,1]$ and a number $m > 0$ such that $f^m(N_{\sigma}) \subset {\rm int} N_{\sigma}$.
\end{cor}
\begin{proofn}
    Due Proposition~\ref{locfl-sep} and Statement~\ref{bicol},  $S_{\sigma}$ is a bi-collared sphere hence there is a topological embedding $e: S^{n-1}\times [-1,1]\to M^n$ such that $e(S^{n-1}\times \{0\})=S_{\sigma}$. Set $N_{\sigma}=e(S^{n-1}\times [-1; 1])$.   
    
    Without loss of generality we suppose that all points of $\partial N_{\sigma}$ belong to the union $W^s_\sigma \cup W^s_\omega$ (otherwise we take $N_{\sigma}$ as the image of $S^{n-1} \times [-\varepsilon,\varepsilon]$  for sufficiently small   $\varepsilon>0$). Then for any point $x\in \partial N_{\sigma}$ there is  $m_x>0$ such that $f^{m_x}(x)\subset {\rm int} N_{\sigma}$. Since $f$ is a homeomorphism,  for any $x\in \partial N_{\sigma}$ there is a neighbourhood $u_x\subset \partial N_{\sigma}$ such that $f^{m_x}(y)\subset {\rm int} N_{\sigma}$ for any $y\in u_x$. Since $\partial N_{\sigma}$ is compact,  the set of neighbourhoods $\{u_x,\, x\in \partial N_{\sigma}\}$ contains a finite subset $\{u_{x_i},\, x_i\in \partial N_{\sigma}, i\in \{1,\dots, \mu\}\}$, covering $\partial N_{\sigma}$.  Set $m=max\{m_{x_i}, i\in \{1,\dots, \mu\}\}$. Then $f^m(\partial N_{\sigma})\subset {\rm int}\, N_{\sigma}$.
    \end{proofn}
\begin{rem}\label{good nbh2} In the case $n=3$ Corollary~\ref{cilnbh} is true but the closure ${\rm cl} W^u_{\sigma^1_2}$ can be a wild sphere in $M^3$, so, the proof of the existence of its neighbourhood is a rather difficult problem. This proof is given in~{\rm \cite{BGMP}} and in~{\rm \cite[Section 6.1.1]{GrPo-book}}. 
\end{rem}

 Suppose that  $\sigma, \omega$ are fixed and  $f(N_{\sigma}) \subset {\rm int}\, N_{\sigma}$ (otherwise  consider the diffeomorphism $f^m$ for an enough big $m\in \mathbb N$). It follows from Lemma~\ref{locfl-sep} and Corollary~\ref{cilnbh}, that the set  $(W^s_\omega\cup W^s_{\sigma})\setminus S_{\sigma}$ consists of two $f$-invariant  connected components $U_+, U_-$. 

\begin{prp} There is a homeomorphism $h_{\pm}: U_\pm\to \mathbb R^n\setminus \{O\} $ such that 
\begin{equation}
    f|_{U_\pm}=h_{\pm}^{-1}a_{0,+1,+1}h_{\pm}|_{U_\pm}.
\end{equation}
\end{prp}
\begin{proofn}
Set $K=N_{\sigma}\setminus {\rm int}\,f(N_{\sigma})$. Since $K$ belongs to an open annulus $S_{\sigma}\times (-1,1)$, it also can be embedded in $\mathbb R^n\setminus \{O\}$. Due to  Annulus theorem (see, for instance~\cite[Theorem 14.3]{GrGu-book_en} for references), $K$ is a union of two disjoint closed annuli $K_+, K_-$.  Suppose that $K_+\subset U_+$. Then  $\bigcup\limits_{n \in \mathbb Z}f^n(K_+)=U_+$ and  for any $x \in (W^s_{\omega} \cup W^s_{\sigma})\setminus S_{\sigma}$ there exists  $n_x \in \mathbb Z$ such that $f^{n_x}(x) \in K_+$.  

Let $S_+=\partial N_\sigma\cap K_+$  and  $\psi_{+}: S_{+} \to \mathbb S^{n-1}$ be an arbitrary homeomorphism.   Define a homeomorphism  $\psi_{1}: f(S_+)\to a_{0,+1,+1}(\mathbb S^{n-1})$ by $\psi_1 = a_{0,+1,+1}\psi_+f^{-1}$.  Then  there exists a  homeomorphism $\psi: K_{+} \to \mathbb K^n$ such that $\psi|_{S_{+}} = \psi_{+}$, $\psi|_{f(S_+)} = \psi_1$ (see ~\cite[Proposition 14.2]{GrGu-book_en} for references). At last, define the  desired homeomorphism $h_+: U_+\to \mathbb R^n\setminus \{O\}$  by 
$h(x)=a_{0,+1,+1}^{-n_x}(\psi(f^{n_x}(x))),$  where  $x\in U_+$   and $f^{n_x}(x)\subset K_+.$  The  homeomorphism $h_-: U_-\to \mathbb R^n\setminus \{O\}$ can be constructed in similar way.
 \end{proofn}

For points $x\in U_\pm, y\in \mathbb R^n\times \mathbb Z_2$ set $x\sim y$ if $y=h_{\pm}(x)$ and denotes by  $M'$  a factor-space of $(M^n\setminus S_{\sigma_{n-1}})\cup (\mathbb R^{n}\times \mathbb Z_2)/_\sim$. The natural projection $p:(M^n\setminus S_{\sigma_{n-1}})\cup (\mathbb R^{n}\times \mathbb Z)\to M'$ induced on $M'$ a structure of a topological  manifold. Denote by $f'$ a map   that coincides with $pf$ on $p(M^n\setminus S_{\sigma})$ and with $pa_{0,+1,+1}$ on each connected component of $p(\mathbb R^n\times \mathbb Z_2)$.  In fact, $M'$ is homeomorphic to a closed manifold, obtained by gluing $M^n\setminus \Int N_\sigma$ and two copies of $\mathbb B^n$ by means of homeomorphisms $h_+, h_-$.  Hence we immediately got the following statement.

\begin{lem}\label{surg-diff}
$M'$ is  homeomorphic to a closed manifold  obtained from $M^n$ by surgery along $S_{\sigma}.$

$f'$ is a regular homeomorphism  of $M'$ and the number $k_i'$ of periodic points of $f'$ having Morse index $i\in \{0,\dots, n\}$   is related to the number $k_i$ of periodic points of $f$ with Morse index $i\in \{0,\dots, n\}$ as follows:

$k'_0=k_0+1, k'_{n-1}=k_{n-1}-1; k'_i=k_i$ for all $i\in \{1,\dots, n-2, n\}$. 
\end{lem}

    We will say that the pair  $\{M',f'\}$ is obtained from $\{M^n, f\}$ by {\it surgery along $W^u_{\sigma}$}.
    
\begin{rem}
    If $M^n$ is a smooth manifold, $f$ is Morse-Smale, and a pair $\{M', f'\}$ is obtained from $\{M^n, f\}$ by the surgery along $W^u_{\sigma_{n-1}}$, then $M'$ is smooth and $f'$ is a Morse-Smale diffeomorphism.
\end{rem}
   
\subsection{End of the proof of Theorem~\ref{main+}}\label{fin}

Let $f\in G(M^n)$ and $k_1^2+k_{n-1}^2\neq 0$. We may suppose that $k_{n-1}\neq 0$ (in the opposite case we consider $f^{-1}$ instead of $f$).  Since  we are interested only in topology of the manifold $M^n$,   we suppose  without loss of generality  that all periodic points of $f$  are fixed (that is 1-periodic, in the opposite case we may consider a homeomorphism $f^N$ for sufficiently large $N$).

Let us remark that if $k_{1}=0$ then similar  to the proof of Proposition~\ref{polar} we obtain that $k_n=1$.

Since $k_{n-1}\neq 0$, there exists a saddle fixed point $\sigma$  of Morse index $(n-1)$ which is the  smallest with respect the Smale relation $\preccurlyeq$ amount of all saddle fixed  points.    Then there is a source $\omega$ such that  $W^u_{\sigma}\setminus \sigma\subset  W^s_\omega$. Due to Lemma~\ref{cilnbh}, the set $S_\sigma=W^u_{\sigma}\cup \{\omega\}$ is a locally flat  sphere. Applying the  surgery operation along $W^u_\sigma$, we obtain a pair $\{f_1, M_1\}$ of a closed topological manifold $M_1$ (may be disconnected) and a regular homeomorphism $f_1$ such that the restriction of $f_1$  on each connected component of $M_1$  belongs to class $G$. If $k_{n-1}=1$  and $k_{1}=0$, then $f_1$ has no saddle  fixed points of indices $1$ and $(n-1)$ and, due to Proposition~\ref{polar}, have only one source. Then  $M_1$ is connected and  $f_1$ is polar.  Due to Lemmas~\ref{surg-diff},~\ref{MeUm-nonor}, $M^n$ is homeomorphic to $M_1\# (S^{n-1}\times S^1)$ if $M^n$ is orientable,  and to to $M_1\# (S^{n-1}\widetilde{\times} S^1)$ if $M^n$ is non-orientable. Since $f_1$ is polar, it has only one sink. Hence, by Lemma~\ref{surg-diff} $k_n=2$. Then $g_{f^t}=1$ and {Theorem~\ref{main+} is proved}.

If $k_{n-1}>1$,  we do the surgery operation until we use up all the saddles of Morse index $(n-1)$ and after $k_{n-1}$  step we obtain a closed manifold $M_{k_{n-1}}$  and a regular homeomorphism $f_{k_{n-1}}: M_{k_{n-1}}\to M_{k_{n-1}}$. There are two possibilities: 1) $k_1=0$; 2) $k_1>0$. In case 1) $k_n=1$ and after each surgery operation we obtain a connected closed  manifold. Then $f_{k_{n-1}}$ is polar, { $M_{k_{n-1}}$} is simply connected and $M^n$ is homeomorphic to $M_{k_{n-1}}\# \underbrace{(S^{n-1}\times S^1) \# \dots \# (S^{n-1}\times S^1)}_{k_{n-1}}$ if $M^n$ is orientable and to  $M_{k_{n-1}}\# \underbrace{(S^{n-1} \widetilde\times S^1) \# \dots \# (S^{n-1}\widetilde{\times} S^1)}_{k_{n-1}}$ otherwise. Since $f_{k_{n-1}}$ is polar, it has only one sink. Hence, by Lemma~\ref{surg-diff} $k_n=1+k_{n-1}$. Then $g_{f^t}=k_{n-1}$ and Theorem~\ref{main+} is proved.

Consider case $2)$.  Continue doing  the surgery operation until we use up all the saddles of Morse index $1$. Then after $\nu=k_{n-1}+k_1$  steps   we obtain a closed manifold $M_{\nu}$  and a regular homeomorphism $f_{\nu}: M_{\nu}\to M_{\nu}$. Let us denote by $\lambda$ the total  number of connected components of $M_\nu$.  Due to Proposition~\ref{polar}, the  restriction of $f_\nu$ on each connected component of $M_\nu$  is polar. Hence the  non-wandering  set of $f_{\nu}$ contains exactly $2\lambda$ sinks and sources. Since the total number of surgery operations is $\nu=k_1+k_{n-1}$, using Proposition~\ref{surg-diff} one obtain that  \begin{equation}\label{2} k_0+k_1+k_{n-1}+k_n=2\lambda.\end{equation} At each surgery operation we have  two possibilities: 1) the operation keeps the number of connected components obtained on the previous steps; 2) the operation increases by one the  number of connected components obtained on the previous steps. Denote by $g$ the number of all operations that have been keeping the number of connected components. Then \begin{equation}\label{3}\lambda=k_1+k_{n-1}-g+1. \end{equation} 

Equations~(\ref{2}),~(\ref{3}) give 
\begin{equation}
g=g_f=(k_1+k_{n-1}-k_0-k_1+2)/2.
\end{equation}
This observation and Lemma~\ref{MeUm-nonor} complete the proof of Theorem~\ref{main+}.

\end{document}